\documentclass[12pt,a4paper]{article}
\usepackage{amsmath}
\usepackage{amssymb}
\usepackage{amsthm}
\usepackage{amsfonts}
\usepackage{latexsym}

\theoremstyle{plain}
\newtheorem{thm}{Theorem}
\newtheorem{cor}{Corollary}
\newtheorem{lem}{Lemma}
\newtheorem{prop}{Proposition}
\theoremstyle{definition}
\newtheorem{defn}{Definition}
\newtheorem{exmp}{Example}

\theoremstyle{remark}
\newtheorem{claim}{Claim}[section]
\newtheorem{rem}{Remark}[section]

\title{Moment maps and equivariant\\Szeg\"o kernels}
\author{Roberto Paoletti\footnote{\noindent{\bf Address.}
Dipartimento di Matematica e Applicazioni, Universit\`a degli
Studi di Milano Bicocca, Via Bicocca degli Arcimboldi 8, 20126
Milano, Italy; {\bf e-mail}: roberto.paoletti@unimib.it }}
\date{}

\begin{document}
\maketitle Let $M$ be a connected $n$-dimensional complex
projective manifold and consider an Hermitian ample holomorphic
line bundle $(L,h_L)$ on $M$. Suppose that the unique compatible
covariant derivative $\nabla _L$ on $L$ has curvature $-2\pi
i\Omega$, where $\Omega$ is a K\"ahler form. Let $G$ be a compact
connected Lie group and $\mu :G\times M\rightarrow M$ a
holomorphic Hamiltonian action on $(M,\Omega)$. Let $\frak{g}$ be
the Lie algebra of $G$, and denote by $\Phi :M\rightarrow
\frak{g}^*$ the moment map.

Let us also assume that the action of $G$ on $M$ linearizes to a
holomorphic action on $L$; given that the action is Hamiltonian,
the obstruction for this is of topological nature \cite{gs-gq}. We
may then also assume that the Hermitian structure $h_L$ of $L$,
and consequently the connection as well, are $G$-invariant.
Therefore for every $k\in \mathbb{N}$ there is an induced linear
representation of $G$ on the space $H^0(M,L^{\otimes k})$ of
global holomorphic sections of $L^{\otimes k}$. This
representation is unitary with respect to the natural Hermitian
structure of $H^0(M,L^{\otimes k})$ (associated to $\Omega$ and
$h_L$ in the standard manner). We may thus decompose
$H^0(M,L^{\otimes k})$ equivariantly according to the irreducible
representations of $G$.

The subject of this paper is the local and global asymptotic
behaviour of certain linear series defined in terms this
decomposition. Namely, we shall first consider the asymptotic
behaviour as $k\rightarrow +\infty$ of the linear subseries of
$H^0(M,L^{\otimes k})$ associated to a single irreducible
representation, and then of the linear subseries associated to a
whole {\it ladder} of irreducible representations. To this end, we
shall estimate the asymptoptic growth, in an appropriate local
sense, of these linear series on some loci in $M$ defined in terms
of the moment map $\Phi$.

To express the problem in point more precisely, fix a Cartan
subalgebra $\frak{h}\subset \frak{g}$, and let $H\subset G$ be the
corresponding Cartan subgroup. Let $R\subset \frak{h}^*$ be the
root system for $(\frak{g},\frak{h})$, and fix a basis $S=\{\alpha
_i\}\subset R$ for $R$. Let $X(H)$ be the group of characters of
$H$, identified with a subgroup of $\frak{h}^*$ (a character $\chi
:H\rightarrow \mathbb{C}^*$ is determined by its differential at
the identity). The finite dimensional irreducible representations
of $G$ correspond bijectively to points in $X(H)\cap
\frak{h}^*_+$, where $\frak{h}_+^*\subset \frak{h}^*$ is the
positive Weyl chamber \cite{serre}. Let $E_\omega$ be the
irreducible representation associated to $\omega \in X(H)\cap
\frak{h}^*_+$. We have for every $k\ge 0$ a $G$-equivariant
orthogonal decomposition
\begin{eqnarray}
H^0(M,L^{\otimes k})\cong \bigoplus _{\omega \in X(H)\cap
\frak{h}^*_+}\, H^0(M,L^{\otimes k})_\omega,
\label{eqn:decomp}\end{eqnarray} where each $H^0(M,L^{\otimes
k})_\omega$ is equivariantly isomorphic to a direct sum of copies
of $E_\omega$. We shall first focus on the asymptotic behaviour,
as $k\rightarrow +\infty$, of the linear series $\left
|H^0(M,L^{\otimes k})_\omega\right |$ for a fixed $\omega \in
X(H)\cap \frak{h}^*_+$, under certain transversality assumptions
on the moment map $\Phi$.

Let us now clarify what is meant here by \textit{local} asymptotic
behaviour of a family of linear series such as $\left
|H^0(M,L^{\otimes k})_\omega\right |$. Suppose first that
$V\subseteq H^0(M,L)$ is any nonzero vector subspace, and let
$\{s_j\}$ be an orthonormal basis of $V$ for the induced unitary
structure. Consider the function
$$\nu _V(p)=\sum _j\left \|s_j(p)\right \|^2_p\,\,\,\,\,\,\,\,\,\,\,\,
(p\in M),$$ where $\|\,\|_p$ denotes the Hermitian norm of $L$ at
$p$. This is in fact independent of the choice of orthonormal
basis, as one can see by replacing the $s_j$'s by $t_k=\sum
_ju_{jk}s_j$ where $U=[u_{ij}]\in U(\dim (V))$ is a unitary
matrix, and then working in a local frame. By the same reason,
given that the representation of $G$ on $H^0(M,L)$ is unitary, if
$V\subseteq H^0(M,L)$ is $G$-invariant then so is the function
$\nu _V$. Here we shall be interested in sequences of functions
such as $\nu _{H^0(M,L^{\otimes k})_\omega}$ and in their
asymptotic behaviour.

Pictorially, we may think of $\nu _V(p)$ as measuring the
\textit{local size}, under the chosen metric structures, of $V$ at
$p$. Asymptotic expansions such as (\ref{eqn:asmypt1}) below have
attracted interest in algebraic geometry since Zelditch remarked
that they control the metric (equivalently, symplectic) asymptotic
behaviour of the maps to projective space associated to the linear
series at hand \cite{z}.

\begin{thm}
Suppose that $0\in \frak{g}^*$ is a regular value of the moment
map $\Phi :M\rightarrow \frak{g}^*$ and that $G$ acts freely on
$\Phi ^{-1}(0)$. Set $n=\dim _{\mathbb{C}}(M)$, $g=\dim
_{\mathbb{R}}(G)$. There exist on $\Phi ^{-1}(0)$ smooth
$G$-invariant real valued functions $a_0>0$ and, for every $\omega
\in X(H)\cap \frak{h}^*_+$ and $l\ge 1$, $a_l^{(\omega)}$ with the
following property: Suppose given $\omega \in X(H)\cap
\frak{h}^*_+$ and the assignment, for every $k\in \mathbb{N}$, of
an orthonormal basis $\left \{s_j^{(k,\omega)}\right \}$ of
$H^0(M,L^{\otimes k})_\omega$. Then there is for $k\gg 0$ an
asymptotic development, uniform in $p\in \Phi ^{-1}(0)$,
\begin{eqnarray} \sum _j\left
|\left |s_j^{(k,\omega)}(p)\right |\right |^2\sim k^{n-g/2}\,\dim
(V_\omega)^2\, a_0(p)\\
+\sum _{l\ge 1} k^{n-l-g/2}a_l^{(\omega)}(p).
\nonumber\label{eqn:asmypt1}
\end{eqnarray}\label{thm:1}
Here $V_\omega$ denotes the irreducible representation associated
to $\omega$.
\end{thm}

If $\tilde G$ denotes the complexification of $G$, the actions of
$G$ on $M$ and $L$ extend to holomorphic actions of $\tilde G$.
Given any subset $A\subseteq M$, we shall denote by $\tilde G\cdot
A$ its saturation under $\tilde G$, that is,
$$\tilde G\cdot A=\left \{g\cdot m:\, g\in \tilde G,\,m\in M\right \}.$$
In the hypothesis of the Theorem, $\tilde G$ acts freely on
$\tilde G\cdot \Phi ^{-1}(0)$, and the latter is the open subset
of stable points for the action of $G$ on $M$ \cite{gs-gq}.

\begin{cor} Under the hypothesis of Theorem \ref{thm:1},
there exists $k_\omega$ such that for $k\ge k_\omega$ the base
locus of the linear series $\left |H^0(M,L^{\otimes k})_\omega
\right |$ satisfies
$${\rm Bs}\left (\left |H^0(M,L^{\otimes k})_\omega \right |\right )
\subseteq M\setminus \left (\tilde G\cdot \Phi ^{-1}(0)\right
).$$\label{cor:1}
\end{cor}

Away from $\Phi ^{-1}(0)$, $\sum _j\left |\left
|s_j^{(k,\omega)}(p)\right |\right |^2$ is rapidly decreasing in
$k$, uniformly so on the complement of the unstable locus. More
precisely, let $R=M\setminus \left (\tilde G\cdot \Phi
^{-1}(0)\right )$ be the set of unstable points of the action.
Then,

\begin{prop}
In the situation of Theorem \ref{thm:1}, if $p\not \in \Phi
^{-1}(0)$ then
$$
\sum _j\left |\left |s_j^{(k,\omega)}(p)\right |\right |^2 =
O(k^{-N}),\,\,\,\,\,N=1,2,\ldots .
$$
This estimate is uniform on compact subsets of $M\setminus \big
(\Phi ^{-1}(0)\cup R\big )$. \label{prop:rapiddecay}
\label{prop:rapiddecay}\end{prop}

Let us at least briefly describe how Theorem \ref{thm:1}
generalizes when the action of $G$ on $\Phi ^{-1}(0)$ is not free.
Since $0$ is a regular value, the action of $G$ on $\Phi ^{-1}(0)$
is at any rate locally free, and therefore the stabilizer subgroup
$G_p$ of any $p\in \Phi ^{-1}(0)$ is finite. Furthermore, there is
an induced unitary action $\alpha _p:G_p\rightarrow S^1$ on the
Hermitian complex line $(L_p,h_p)$. We may then take the
$L^2$-product of $\alpha _p$ and $\chi _\omega$ on $G_p$ with
respect to the counting measure, that is,
$$(\chi _\omega,\alpha _p)_{G_p}=\sum _{g\in G_p}\chi _\omega
(g)\cdot \overline {\alpha _p(g)}.$$ Then (\ref{eqn:asmypt1}) is
replaced by
$$\sum _j\left
|\left |s_j^{(k,\omega)}(p)\right |\right |^2\sim k^{n-g/2}\,\dim
(V_\omega)\,(\chi _\omega,\alpha _p^k)_{G_p}\, a_0(p)+\sum _{l\ge
1} k^{n-l-g/2}a_l^{(\omega)}(p),$$ uniformly in $p\in \Phi
^{-1}(0)$.

The microlocal techniques used in the proof of Theorem \ref{thm:1}
can also be used to study the asymptotic growth of the dimension
of the spaces of equivariant sections $H^0(M,L^{\otimes k})
_\omega$ when $\omega$ is kept fixed and $k\rightarrow +\infty$.
There already exist two approaches to this problem, one geometric
and the other algebraic.

A geometric solution follows from Meinrenken's proof of a
fundamental conjecture of Guillemin and Sternberg \cite{mein},
resting on the symplectic cutting technique of E. Lerman
\cite{ler}: Obviously $G$ acts on every cohomology group
$H^i(M,L)$, $i=0,\ldots,n$, and we may consider the virtual vector
space
$$ \mathrm{RR}(M,L)_\omega=\sum _{i=0}^n(-1)^iH^i(M,L)_\omega.
$$ Let $\mu _\omega =\mu _\omega (M,L)$ be the corresponding virtual
multiplicity for the representation $V_\omega$. By the main
Theorem of \cite{mein}, $\mu _\omega$ may be computed as a
Riemann-Roch number on the symplectic reduction $M_\omega$ of the
Hamiltonian $G$-manifold $(M,\Omega,\Phi)$ at $\omega$, provided
$\omega$ is a regular value of the moment map. Then $M_\omega
=\Phi ^{-1}(\omega)/G_\omega$, where $G_\omega\subseteq G$ is the
stabilizer subgroup of $\omega$ under the coadjoint action.
Assuming to fix ideas that $G_f$ acts freely on $\Phi
^{-1}(\omega)$, $M_\omega$ naturally inherits by restriction and
quotient a K\"{a}hler structure and a compatible polarization
$L_\omega$. Then by Theorem 1.1 and Corollary 1.2 of \cite{mein}
$$\mu _\omega =\mathrm{RR}(M_\omega,L_\omega).$$
Let us now replace $L$ by $L^{\otimes k}$, and thus $\Omega$ by
$k\Omega$ and $\Phi$ by $\Phi _k=k\Phi$. By Serre vanishing if
$k\gg 0$ then $H^i(M,L^{\otimes k})=0$ for all $i>0$, and
therefore $ \mathrm{RR}(M,L^{\otimes k})_\omega=H^0(M,L^{\otimes k
})_\omega$ and $\mu _\omega (L^{\otimes k})$ is the multiplicity
of $V_\omega$ in $H^0(M,L^{\otimes k })$. Thus,
$$\dim H^0(M,L^{\otimes k })_\omega=\dim (V_\omega)\cdot
\mathrm{RR}(M_{\omega ,k},L_{\omega ,k}^{\otimes k}),$$ where
$M_{\omega ,k}=\Phi ^{-1}(k^{-1}\omega)/G_\omega$ is the reduction
of $(M,k\Omega,\Phi _k)$ at $\omega$ and $L_{\omega ,k}^{\otimes
k}$ is the polarization on it induced by $L^{\otimes k}$. Now by
Corollary 7.3 of \cite{mein} the reduced spaces $M_{\omega ,k}$
are all diffeomorphic to the fibration
$$M(P,\mathcal{O}_\omega)=P\times _G\mathcal{O}_\omega,$$ where
$P=\Phi ^{-1}(0)$ is viewed as a principal $G$-bundle over the
symplectic quotient $M_0=\Phi ^{-1}(0)/G$, and
$\mathcal{O}_\omega$ is the coadjoint orbit of $\omega$. On the
other hand, again by Corollary 7.3 of \cite{mein}, the symplectic
structure $\Omega _{\omega ,k}$ on $M_{\omega ,k}$ is the one
induced by minimal coupling \cite{gs} from $k\Omega$, the
symplectic structure $\sigma _{k^{-1}\omega}$ on
$\mathcal{O}_{k^{-1}\omega}$, a fixed connection on $P$ and the
moment map $J_{\omega, k}=\frac
1kJ_\omega:\mathcal{O}_\omega\hookrightarrow \frak{g}^*$ of
$\mathcal{O}_{\omega,k}$. With the appropriate scaling taken into
account, one obtains a leading asymptotics of degree $n-g$.

There is also an algebraic line of research on this asymptotic
problem, in the work of Brion and Dixmier \cite{bd}, \cite{b}.

Here we propose a different, analytic and fairly elementary
approach to the same asymptotics: namely, we use Boutet de Monvel
and Sj\"{o}strand's microlocal description of the Szeg\"{o} kernel
\cite{bs} to reduce the problem to an application of the
stationary phase Lemma along $\Phi ^{-1}(0)$. Behind some
technicalities, the basic idea is very simple.

In this formulation, the asymptotics depends on the weight
$\omega$ only through the germ at the identity $e\in G$ of its
charachter function $\chi _\omega \in \mathcal{C}^\infty (G)$,
rather than on the geometry of the coadjoint orbit
$\mathcal{O}_\omega$.

\begin{thm} In the hypothesis and notation of Theorem \ref{thm:1},
there exist for $j\ge 1$ differential polynomials $S_j$ of degree
$2j$ on a neighbourhood of $e\in G$ such that for any $\omega \in
X(H)\cap \frak{h}^*_+$ we have
\begin{eqnarray} \dim H^0(M,L^{\otimes k}) _\omega \,\thicksim \,
\dim (V_\omega)^2\cdot \mathrm{vol}(M_0) \cdot k^{n-g}\\
+\dim (V_\omega)\,\sum _{j\ge 1}b_{j,\omega}\,k^{n-g-j},\nonumber
\label{eqn:asymptoticsections}
\end{eqnarray} where
$\mathrm{vol}(M_0)$ is the volume of the reduced space $M_0=:\Phi
^{-1}(0)/G$ with its natural K\"{a}hler structure, and
$b_{j,\omega}=S_j(\chi _\omega)(e)$. \label{thm:3}\end{thm}

It follows from the proof of Theorem \ref{thm:3} (and the
stationary phase Lemma) that the $S_j$'s may be expressed in terms
of the Hessian of $h_L$ along $\Phi ^{-1}(0)$ and the classical
symbol appearing in the Fourier integral representation of the
Szeg\"{o} kernel of $\Pi$ given in \cite{bs}.

This extends to our setting a (special case of a) result of Brion
and Dixmier \cite{bd}, \cite{b}:

\begin{cor} In the hypothesis of
Theorem \ref{thm:1}, let $\mu _{\omega ,k}$ be the multiplicity of
$V_\omega$ in $H^0(M,L^{\otimes k})_\omega$. Then
$$\lim _{k\rightarrow +\infty}\frac{\mu _{\omega ,k}}{\mu _{0 ,k }} =
\dim (V_\omega).$$\end{cor}

When $0\not\in \Phi (M)$, a more informative result is provided by
the study of the corresponding asymptotic properties of the {\it
ladder} linear series
$$H^0(M,L^{\otimes k})_{(\omega)}=:
\bigoplus _{\ell=1}^\infty H^0(M,L^{\otimes k})_{\ell
\omega}\,\subseteq \,H^0(M,L^{\otimes k}).$$ We shall now assume
that $G$ is semisimple.

The action of $G$ on $L$ induces an action on the dual line bundle
$L^*$, equipped with the dual hermitian metric. The unit circle
bundle $X\subseteq L^*$ is invariant under $G$, and therefore
there is an induced Hamiltonian action of $G$ on the cotangent
bundle $T^*X$, and on the complement of the zero section
$T^*X\setminus \{0\}$. Let $\Psi :T^*X\setminus \{0\}\rightarrow
\frak{g}^*$ be the corresponding (conic) moment map.

Since the connection is $G$-invariant, the action on $T^*X$
preserves the positive cone $Y\subseteq T^*X$ generated by the
normalized connection 1-form $\alpha$:
\begin{equation}Y=\{(x,r\alpha _x):x\in
X,r>0\}.\label{eqn:Y}\end{equation} It is well-known that since
$\Omega$ is symplectic $Y$ is in fact a symplectic conic
submanifold of $T^*X$.

The moment maps $\Psi$ and $\Phi$ are tightly related on $Y$: one
has \begin{equation}\Psi \left ((x,r\alpha _x)\right )=r\Phi (\pi
(x)),\,\,\,\,\,\,\,\,\,\,(x\in X)\label{eqn:momY} \end{equation}
where $\pi :X\rightarrow M$ is the projection \cite{gs-gq}.

We need a further piece of notation \cite{gs-hq}.

\begin{defn} For $\omega \in X(H)\cap \frak{h}^*_+$, we shall denote by
$\mathbb{R}_+\omega =\{r\cdot \omega: r>0\}$ the positive ray
through $\omega$. Let furthermore ${\cal O} ={\cal O}_\omega
\subseteq \frak{g}^*$ be the coadjoint orbit of $\omega$. We shall
denote by $C({\cal O})\subseteq \frak{g}^*$ the positive cone over
${\cal O}$, that is,
$$C({\cal O})=\left \{r \lambda:\, r>0,\, \lambda \in {\cal O}\right \}.$$
\label{defn:cone}\end{defn}

Notice that $\Phi :M\rightarrow \frak{g}^*$ is transversal to
$C(\mathcal{O})$ if and only if so is $\left .\Psi \right |_Y$.
Since $\left .\Psi \right |_Y$ is conic, the latter condition is
in turn equivalent to the one that $\left .\Psi \right |_Y$ be
transversal to $\mathcal{O}$, and by $G$-equivariance this is in
turn equivalent to the condition that $\omega$ be a regular value
of $\left .\Psi \right |_Y$. Hence, in view of the invariance of
the connection 1-form, conditions ii) and iii) below are
equivalent to the one that $\omega$ should lie in an
\textit{elementary fundamental wedge} for the induced Hamiltonian
action of $G$ on $Y$, in the terminology of \cite{gs-hq}, page
357.

\begin{thm}\label{thm:2}
Suppose $G$ is semisimple, compact and connected. Given $\omega
\in X(H)\cap \frak{h}^*_+$, suppose that

\noindent i): $0\not\in \Phi (M)$;

\noindent ii): $\Phi$ is transversal to $C({\cal O})$;

\noindent iii): the stabilizer subgroup $G_\omega\subseteq G$ of
$\omega$ acts freely on $\Phi ^{-1}\big (\mathbb{R}_+\omega\big )$
(equivalently, $G_\omega$ acts freely on $\Phi ^{-1}\big (C({\cal
O})\big )$).

\noindent Then there exist smooth $G$-invariant functions
$$b_l:\Phi ^{-1}\left (C({\cal O})\right ) \rightarrow \mathbb{R}
\,\,\,\,\,\,\,\,\,\,\,\,
\text{\textrm{($l=0,1,\ldots$)}}$$ with $b_0>0$ such that for
every choice of orthonormal basis $\left \{s_j^{(k,\ell
\omega)}\right \}$ of $H^0(M,L^{\otimes k})_{(\ell \omega)}$ for
$\ell =1,2,\ldots$, there is an asymptotic development
\begin{eqnarray}
\sum _{\ell,j}\left |\left |s_j^{(k,\ell \omega)}(p)\right |\right
|^2\sim \sum _{l\ge 0}k^{n-l}\,b_l(p),
\end{eqnarray}
uniformly in $p\in \Phi ^{-1}\left (C({\cal O})\right )$.
Furthermore,
\begin{eqnarray}
\sum _{\ell,j}\left |\left |s_j^{(k,\ell \omega)}(p)\right |\right
|^2= O(k^{-N}),\,\,\,\,\,N=1,2,\ldots,
\end{eqnarray}
uniformly on compact subsets of $M\setminus \Phi ^{-1}\left
(C({\cal O})\right )$. \end{thm}

\begin{cor} Under the hypothesis of Theorem \ref{thm:2},
there exists $k_\omega$ such that for $k\ge k_\omega$ the base
locus of the linear series $\left |H^0(M,L^{\otimes k})_{(\omega)}
\right |$ satisfies
$${\rm Bs}\left (\left |H^0(M,L^{\otimes k})_{(\omega)}
\right |\right )\subseteq M\setminus \left (\tilde G\cdot W\right
),$$ where $W=\Phi ^{-1}\left (C({\cal O})\right )$.\label{cor:2}
\end{cor}

By restriction of the arguments in the proof of Theorem we have:

\begin{cor} If the hypothesis of Theorem \ref{thm:2} are satisfied on a
$G$-invariant open set $M'=M\setminus B$, where $B\subseteq M$ has
measure zero, then the same conclusions hold for $p\in
M'$.\label{cor:M'}
\end{cor}

\begin{exmp} If $\Omega _{FS}$ is the Fubini-Study form on
$\mathbb{P}^1$, let $\Omega =: 2\pi _1^*\left (\Omega _{FS}\right
)+\pi _2^*\left ( \Omega _{FS}\right )$, where $\pi
_i:\mathbb{P}^1\times \mathbb{P}^1\rightarrow \mathbb{P}^1$ is the
projection on the $i$-th factor ($i=1,2$). Consider the diagonal
action of $\mathrm{SU}(2)$ on $(\mathbb{P}^1\times \mathbb{P}^1,
\Omega)$. The moment map $\hat{\Phi}:\mathbb{P}^1\times
\mathbb{P}^1\rightarrow \frak{su}(2)^*$ for this action is
$\hat{\Phi}(p,q)=2\Phi (p)+\Phi (q)$, where $\Phi
:\mathbb{P}^1\rightarrow \frak{su}(2)^*$ is the moment map for the
action on $(\mathbb{P}^1, \Omega _{FS})$. We may equivariantly
identify $\mathbb{P}^1$ with the unit sphere $S^2$ and
$\frak{su}(2)^*\cong \frak{su}(2)\cong \mathbb{R}^3$, so that
$\Phi$ corresponds to the inclusion $\iota :S^2\hookrightarrow
\mathbb{R}^3$. Then the hypothesis in Theorem \ref{thm:2} and its
Corollary \ref{cor:M'} are satisfied on the complement of the
diagonal, $M'=\left (\mathbb{P}^1\times \mathbb{P}^1\right )
\setminus \Delta _{\mathbb{P}^1}$ (strictly speaking, after
replacing the action of $\mathrm{SU}(2)$ on $\mathbb{P}^1$ by the
action of $\mathrm{SO}(3)$ on $S^2$). If $H$ denotes the
hyperplane line bundle on $\mathbb{P}^1$, the line bundle on
$\mathbb{P}^1\times \mathbb{P}^1$ associated to $\Omega$ is
$L=H^{\otimes 2}\boxtimes H$. Set $V=\mathbb{C}^2$. For every
$k\ge 1$,
\begin{eqnarray*}H^0(\mathbb{P}^1\times \mathbb{P}^1, L^{\otimes k})\cong
\mathrm{Sym}^{2k}(V^*)\otimes \mathrm{Sym}^k(V^*) \cong \bigoplus
_{0\le j\le k}\mathrm{Sym}^{3k-2j}(V^*).
\end{eqnarray*}
For any integer $r\ge 1$, let $H^0(\mathbb{P}^1\times
\mathbb{P}^1, L^{\otimes k})_{(r)}\subseteq H^0(\mathbb{P}^1\times
\mathbb{P}^1, L^{\otimes k})$ be the subspace corresponding to the
direct sum of the terms $\mathrm{Sym}^{3k-2j}(V^*)$ with
$3k-2j\equiv 0$ ($\mathrm{mod}\, r$). By the Theorem and its
Corollary, for any fixed $(p,q)\in \left (\mathbb{P}^1\times
\mathbb{P}^1\right )\setminus \Delta _{\mathbb{P}^1}$ and any
fixed integer $r\ge 1$ the asymptotic growth as $k\rightarrow
+\infty$ of the local size at $(p,q)$ of the linear series
$H^0(\mathbb{P}^1\times \mathbb{P}^1, L^{\otimes k})_{(r)}$ grows
like $a_r(p,q)k^2$ for some $a_r(p,q)>0$.
\end{exmp}

\bigskip

These results are based on the microlocal description of the
Szeg\"o kernel given by Boutet de Monvel and Sj\"ostrand in
\cite{bs}. In particular, Theorem \ref{thm:2} is proved by giving
a similar microlocal description of the orthogonal projector
associated to the linear series $H^0(M,L^{\otimes k})_{(\ell
\omega)}$. To this end, we shall also rely on the microlocal
description of a projector associated to ladders of
representations given by Guillemin and Sternberg in \cite{gs-hq},
and on a reduction technique used by Schiffman and Zelditch in
\cite{sz}. Overall, as the reader will easily see, the paper is
also largely in debt to arguments from \cite{gs-gq}, \cite{sz} and
\cite{z}.

In future work we shall consider extensions of these results to
the almost complex setting, and further investigate the asymptotic
growth of the spaces of equivariant sections.

\bigskip

\noindent \textbf{Acknowledgments.} I am endebted to the referee
for some fruitful and stimulating comments.

\section{Proof of Theorem \ref{thm:1}}
\noindent Let $L^*$ be the dual line bundle of $L$, with the
induced Hermitian metric and connection, and let $$\rho
:L^*\rightarrow \mathbb{R},\,\,\,\,\,(x,v)\mapsto ||v||_x^2,$$ be
the associated square norm function. Let $X\subset L^*$ be the
unit circle bundle:
$$X=\left \{(p,v)\,:\,\rho \left ((p,v)\right )=1\right \},$$
with projection $\pi :X\rightarrow M$. We shall denote by
$p,q,\ldots$ points in $M$, and by $x,y,\ldots$ points in $X$.

By the ampleness of $L$, $X$ is the boundary of the bounded
strictly pseudoconvex domain $D =\{\rho \le 1\}$. If $i\, \alpha$
is the connection form on $X$, then $\alpha$ is a contact form,
$d\alpha =\pi ^*(\Omega)$ and $\alpha \wedge \pi
^*(\Omega)^{\wedge n}$ is a volume form on $X$. Given this, we
shall implicitly identify functions and half-forms.

There is a canonical isomorphism for every $k$ between the spaces
of smooth sections of $L^{\otimes k}$ on $M$, $\mathcal{C}^\infty
(M,L^{\otimes k})$, and the spaces $\mathcal{C}^\infty (X)_k$ of
smooth functions on $X$ of the $k$-th isotype for the
$S^1$-action. We shall occasionally denote by $\tilde V \subseteq
\mathcal{C}^\infty (X)_k $ the subspace corresponding to a
subspace $V\subseteq \mathcal{C}^\infty (M,L^{\otimes k})$, and
occasionally not distinguish between the two.

Let $\Pi \in {\cal D}'(X\times X)$ be the Szeg\"o kernel, that is,
the distributional kernel of the orthogonal projector $\pi
:L^2(X)\rightarrow H(X)$, where $H(X)$ is the Hardy space of
boundary values on $X$ of holomorphic functions on $D$. If
$\{s_j^{(k)}\}_{j=1}^{N_k}$ is an orthonormal basis of
$H^0(M,L^{\otimes k})$ for every $k\ge 0$, viewed implicitly as a
space of CR functions on $X$, we have
$$\Pi (x,y)=\sum _{k=0}^{+\infty}\, \Pi _k(x,y)\,\,\,\,\,\,(x,y\in X),$$ where
\begin{equation}\Pi _k(x,y)=\sum _{j=1}^{N_k}\,
s_j^{(k)}(x)\otimes \overline s_j^{(k)}(y)\,\,\,\,\,\,\,\,
\mbox{($x,y\in X,\,k\ge 0$)}.\label{eqn:pik}\end{equation} As
proved in \cite{bs}, $\Pi$ is a Fourier integral with complex
phase. More precisely, it is microlocally equivalent to an
oscillatory integral of the form \begin{equation}\Pi (x,y)=\int
_0^{+\infty}\, e^{i\, t\, \psi (x,y)}\,
s(x,y,t)\,dt\,\,\,\,\,\,(x,y\in
X),\label{eqn:fourierintegral}\end{equation} where $s\in
S^n(X\times X\times \mathbb{R}_+)$ has an asymptotic expansion
\begin{equation}
s(x,y,t)\sim \sum _{j=0}^\infty s_j(x,y)\, t^{n-j}.
\label{eqn:asymptotic}\end{equation} The restriction to the
diagonal of the principal term $s_0(x,x)$ is given explicitly in
equation (4.10) of \cite{bs}:
\begin{equation}s_0(x,x)=\frac 1{4\pi ^n}\det \left (L_X(x)\right )
\cdot ||d\rho||\,\,\,\,\,(x\in X),\label{eqn:s0}
\end{equation}
where $L_X$ is the Levi form. We refer to \cite{bs}, \cite{z},
\cite{sz} for a discussion of the phase $\psi \in {\cal C}^\infty
(L^*\times L^*)$; it parametrizes an almost holomorphic Lagrangian
submanifold, whose real locus is the wave front of $\Pi$. This is
the isotropic conic submanifold
\begin{eqnarray}\Sigma =\left \{(x,r\alpha _x,x,-r\alpha _x):
\, x\in X,\, r>0\right \}\subseteq T^*(X\times
X).\label{eqn:wf}\end{eqnarray} The Taylor series of $\psi$ along
the diagonal $\Delta _{L^*} \subset L^*\times L^*$ is completely
determined (equivalently, $\psi$ is uniquely determined up to a
function vanishing to infinite order along $\Delta _{L^*}$).
Explicitly, if $x\in L^*$, in local holomorphic coordinates
induced by a local holomorphic frame for $L$ in a neighbourhood of
$\pi (x)$, we have
\begin{eqnarray}\psi (x+h,x+k)\, \sim \, \frac 1i\,\sum _{I,J}\,
\frac{\partial ^{I+J}\rho}{\partial z^I\partial \overline z^J}(x)
\,\frac{h^I}{I!}\, \frac{\overline k^J}{J!}
\,\,\,\,\,\,\,\,\,\,\,(h,k\in
\mathbb{C}^{n+1}).\label{eqn:psi}\end{eqnarray} We can retrieve
$\Pi _k$ in (\ref{eqn:pik}) as the $k$-th Fourier component of
$\Pi$:
$$\Pi _k(x,y)=\int _0^{+\infty}\,\int _{S^1}e^{-ik\theta}
e^{it\psi (r_\theta x,y)}s(r_\theta x,y,t)\,dt\,d\theta
\,\,\,\,\,(x,y\in X),$$ where $r:(e^{i\theta},x)\in S^1\times
X\mapsto r_\theta (x)\in X$ is the $S^1$-action on $X$;
application of the stationary phase lemma gives an asymptotic
expansion for $\Pi _k$ in terms of which many classical results in
algebraic geometry can be deduced \cite{z}.

Given the direct sum decomposition (\ref{eqn:decomp}), we may take
as an orthonormal basis for $H^0(M,L^{\otimes k})$ the union of a
collection of orthonormal basis
$\{s_j^{(k,\omega)}\}_{j=1}^{N_{k,\omega}}$ of $H^0(M,L^{\otimes
k})_\omega$ for each $\omega \in X(H)\cap \frak{h}^*_+$. Thus,
$$\Pi _k(x,y)=\sum _{\omega \in X(H)\cap \frak{h}^*_+}\,
\Pi _{k,\omega}(x,y)\,\,\,\,\,\,\,\,\mbox{($x,y\in X,\, k\ge
0$)},$$ where $$\Pi _{k,\omega}(x,y)=\sum _{j=1}^{N_{k,\omega}}
\,s_j^{(k,\omega)}(x)\otimes \overline s_j^{(k,\omega)}(y)
\,\,\,\,\,\,\,\,\mbox{($k\ge 0,\,\omega \in X(H)\cap
\frak{h}^*_+$)}.$$ If $x\in X$ and $p=\pi (x)$, we have
$$\Pi _{k,\omega}(x,x)=\sum _{j=1}^{N_{k,\omega}}
\,||s_j^{(k,\omega)}(p)||^2.$$ Thus we want to study the
asymptotic behaviour of $\Pi _{k,\omega}(x,x)$ for a fixed
$\omega$ and $x\in \pi ^{-1}\left (\Phi ^{-1}(0)\right )$ as
$k\rightarrow +\infty$.

Clearly, $\Pi _{k,\omega}$ is the $k$-th Fourier component of the
equivariant Szeg\"{o} kernel $\Pi _\omega \in \mathcal{D}'(X\times
X)$ associated to $\omega$, that is, the distributional kernel of
the orthogonal projector
$$\pi _\omega :L^2(X)\longrightarrow H(X)_\omega,$$
where we have set $$H(X)_\omega = \bigoplus _{k\in
\mathbb{Z}}\widetilde{H^0\left (M,L^{\otimes k }\right
)_\omega}.$$

Thus, $\Pi _{\omega ,k}\in \mathcal{C}^\infty (X\times X)$ is the
distributional kernel for the orthogonal projector $$\pi _{\omega
,k }:L^2(X)\longrightarrow \widetilde{H^0\left (M,L^{\otimes
k}\right )_\omega}.$$

Let $p_\omega :L^2(X)\rightarrow L^2(X)_\omega$ be the orthogonal
projector onto the $G$-equivariant Hilbert subspace of $L^2(X)$
associated to $\omega$, let $q_k:L^2(X)\rightarrow L^2(X)_k$ be
the orthogonal projector onto the $k$-th isotype for the
$S^1$-action, and let $\pi :L^2(X)\rightarrow H(X)$ be as above
the orthogonal projector onto the Hardy space of CR functions.
Then
\begin{eqnarray}\pi _{\omega ,k}\, =\, p_\omega \circ q_k\circ
\pi.\label{eqn:compos}\end{eqnarray}

Let $G$ be a compact topological group and $\sigma :G\rightarrow
U(V)$ be a finite dimensional irreducible representation. Let
$\varrho :G\rightarrow U(W)$ be a unitary action on a separable
Hilbert space, and let $W_V\subseteq W$ be the equivariant
$G$-suspace associated to $\sigma$. Then the orthogonal projection
operator $\Pi _V:W\rightarrow W_V$ is given by
$$\Pi _V=\dim (V)\cdot \int _G \, \varrho (g)\, \chi _\sigma (g^{-1})\, dg,$$
where $\chi _\sigma$ is the character of $\sigma$ and $dg$ a
normalized Haar measure on $G$ \cite{dixmier}.

To simplify notation, let us denote by $\tilde \mu $ both the
linearization to $L$ of the action of $G$ on $M$, and the induced
actions of $G$ on $L^*$ and $X$. Under the usual identifications,
the linear representations of $G$ on the spaces $H^0(M,L^{\otimes
k})$ are then given by pull-back of CR functions on $X$ under this
action. Thus,
\begin{eqnarray}\Pi _{k,\omega}(x,y)=
\dim (V_\omega)\cdot \int _G \,  \chi _\omega (g^{-1})\,
\varrho (g)\left (\Pi _k(x,y)\right )\, dg \nonumber \\
=\dim (V_\omega)\cdot \int _G \,  \chi _\omega (g^{-1})\,
\Pi _k(\tilde \mu _{g^{-1}}(x),y)\, dg \nonumber \\
=\dim (V_\omega)\cdot\int _0^{+\infty}\,\int _{S^1}\,\int _G \,
e^{-ik\theta}e^{it\psi \big (\tilde \mu _{g^{-1}}\circ r_\theta
(x),y\big )} \cdot  \chi _\omega (g^{-1})\cdot \nonumber \\ \cdot
s\big (\tilde \mu _{g^{-1}}\circ r_\theta (x),y,t \big
)\,dt\,d\theta \,dg\,\,\,\,\,\,\,\,\,\,\,\,\,(x,y\in X)
\label{eqn:komega}\end{eqnarray} By assumption, $G$ acts freely on
$\Phi ^{-1}(0)$. We may thus find an open neighbourhood $U$ of the
unit $e\in G$ and $\epsilon >0$ such that $d(p,gp)>\epsilon$ if
$p\in \Phi ^{-1}(0)$, $g\not\in U$. If $x=y\in \pi ^{-1}\left
(\Phi ^{-1}(0)\right )$, the integral (\ref{eqn:komega}) may thus
be decomposed as follows:
\begin{eqnarray}\Pi _{k,\omega}(x,x)=\dim (V_\omega )\cdot \int _U \,
\chi _\sigma (g^{-1})\,\Pi _k(\tilde \mu _{g^{-1}}(x),x)\, dg \nonumber \\
+\dim (V_\omega)\cdot \int _{G\setminus U} \,  \chi _\sigma
(g^{-1})\, \Pi _k(\tilde \mu _{g^{-1}}(x),x)\, dg,
\label{eqn:komega1}\end{eqnarray} the latter term being
$O(k^{-N})$ for every $N=1,2,\ldots$ as $k\rightarrow +\infty$,
uniformly in $x\in \pi ^{-1}\left (\Phi ^{-1}(0)\right )$. Let us
focus on the former term, which we call $\Pi _{k,\omega}(x,x)'$.
Setting $t=ku$, this may be rewritten
\begin{eqnarray}\Pi _{k,\omega}(x,x)'=k\,\dim (V_\omega )\cdot\int _0^{+\infty}\,
\int _{S^1}\,\int _U \, e^{ik\big [u\psi \big (\tilde \mu
_{g^{-1}}\circ r_\theta (x),x\big )
-\theta \big ]}\cdot   \chi _\sigma (g^{-1})\cdot \nonumber \\
\cdot\,s\big (\tilde \mu _{g^{-1}}\circ r_\theta (x),x,ku\big )
\,du\,d\theta \,dg \nonumber \\= k\,\dim (V_\omega)\cdot\int
_0^{+\infty}\int _{S^1}\int _U  e^{ik\Psi (x,u,\theta,g)} \cdot
\chi _\sigma (g^{-1})\cdot \nonumber \\
\cdot s \big (\tilde \mu _{g^{-1}}\circ r_\theta (x),x,ku\big
)\,du\,d\theta \,dg. \label{eqn:komega2}\end{eqnarray} From the
corresponding property of $\psi$, it follows that the phase
$$\Psi (x,u,\theta,g)=u\psi \big (\tilde \mu _{g^{-1}}\circ r_\theta (x),x\big )
-\theta$$ has positive imaginary part. Therefore,
(\ref{eqn:komega2}) is a complex oscillatory integral, and its
asymptotic behaviour as $k\rightarrow +\infty$ is determined by
the stationary points of the phase  as a function of
$(g,\theta,t)$.

\begin{lem}\label{lem:critical}
Suppose $x\in \Phi ^{-1}(0)$. Then $(e,0,1)$ is a non-degenerate
critical point of $\Psi$. Furthermore, perhaps after replacing $U$
with a smaller open neighbourhood of $e\in G$, it is the only
critical point of $\Psi$ in $U\times S^1\times
(0,+\infty)$.\end{lem}

Here, of course, we implicitly identify $\theta$ with
$e^{i\theta}$.

\bigskip

\noindent {\it Proof of Lemma \ref{lem:critical}.} Let us first
show that $(e,0,1)$ is a critical point of $\Psi$. We have
\begin{equation}\label{eqn:initial} \tilde \mu _e\circ
r_0(x)=x\mbox{ and }(d\psi )_{(x,x)}=(x,\alpha _x,x,-\alpha
_x).\end{equation} The connection on $L$ induces a $G\times
S^1$-invariant direct sum decomposition $TX=H(X/M)\oplus V(X/M)$
into a horizontal and a vertical subbundle; here
\begin{center}$H(X/M)=\ker (\alpha)$ and $V(X/M)={\rm span}\left
\{\frac \partial {\partial \theta}\right \}$\end{center} ($\frac
\partial {\partial \theta}$ denotes the generator of the
$S^1$-action on $X$). If $\xi \in \frak{g}$, let $\xi _M$ and $\xi
_X$ denote, respectively, the vector fields induced by $\xi$ on
$M$ and $X$. If $V$ is a vector field on $M$, let $V^\sharp$
denote its horizontal lift to $X$. Then, in terms of the above
direct sum decomposition of $TX$, we have for all $\xi \in
\frak{g}$ \cite{gs-gq}:
\begin{eqnarray}
\xi _X=\left (\xi _M^\sharp,\, (\phi _\xi \circ \pi) \cdot \frac
\partial {\partial \theta}\right )\label{eqn:hor}
\end{eqnarray}
where $\phi _\xi =:\left <\Phi,\,\xi \,\right >:M\rightarrow
\mathbb{R}$ is the $\xi$-component of the moment map. In
particular, \begin{equation}\xi _X(x)=\xi _M^\sharp (x)\mbox{ for
every $\xi\in \frak{g}$ if }\Phi (\pi (x))=0.
\label{eqn:phi0}\end{equation} Let us introduce the map
$$a_\xi :\mathbb{R}\rightarrow X\times X,\,\,\,\,\,\,t\mapsto
\left (\tilde \mu _{\exp(t\xi)}(x),x\right ).$$ Given
(\ref{eqn:initial}) and (\ref{eqn:phi0}),
\begin{eqnarray}(\partial _\xi \Psi)_{(e,1,0)}=\left .
\frac d{dt}\right |_{t=0}\psi \left (a_\xi (-t)x,x\right ) =
(d\psi )_{(x,x)}\left (-\xi _X(x),0\right )=0,\end{eqnarray} for
every $\xi \in \frak{g}$. On the other hand, it follows from the
arguments on pages 327-328 of \cite{z} that $(1,0)$ is the only
critical point of $\Psi (x,e,\theta,u)$ as a function of
$(\theta,u)$; more precisely, we have
\begin{eqnarray}\label{eqn:critpsi}
\left (d_u\Psi \right )_{(x,e,\theta,u)}=-i\,(1-e^{i\theta})\mbox{
and } \left (d_\theta\Psi \right
)_{(x,e,\theta,u)}=te^{i\theta}-1.\end{eqnarray} Furthermore
$(1,0)$ is a nondegenerate critical point of $\Psi (x,e,\theta,u)$
as a function of $(\theta,u)$, and the Hessian is given there by
\begin{equation}\label{eqn:hessian}\left [\begin{array}{cc}i&1\\1&0\end{array}\right ].
\end{equation}
We can now prove that $(e,0,1)$ is a nondegenerate critical point
of $\Psi$. Perhaps after restricting to a smaller open
neighbourhood of $e\in G$, we may suppose that $U$ is
diffeomorphic to an open neighbourhood of $0\in \frak{g}$ under
the exponential map; having fixed a basis of $\frak{g}$, let
$(h_1,\ldots,h_g)$ ($g=\dim (G)$) be the resulting coordinates on
$U$ centered at $e$. By the above and the expression of $\Psi$, it
is clear that $\left .\frac {\partial ^2\Psi}{\partial h_j\partial
u}\right |_{(e,1,0)}=0$. Given this and (\ref{eqn:hessian}), the
Hessian of $\Psi (x,g,\theta,u)$ in the coordinates
$(h_j,\theta,u)$ at $(e,0,1)$ has the form
$$H(\Psi)=\left [\begin{array}{ccccc}\left .\frac {\partial ^2\Psi}{\partial h_1^2}
\right |_{(e,0,1)}&\cdots&\left .\frac {\partial ^2\Psi}{\partial
h_1\partial h_g} \right |_{(e,0,1)}&\left .\frac {\partial
^2\Psi}{\partial h_1\partial \theta}
\right |_{(e,0,1)}&0\\ \vdots&\ddots&\vdots&\ddots&\vdots\\
\left .\frac {\partial ^2\Psi}{\partial h_g\partial h_1}\right
|_{(e,0,1)}&\cdots& \left .\frac {\partial ^2\Psi}{\partial
h_g^2}\right |_{(e,0,1)}&
\left .\frac {\partial ^2\Psi}{\partial h_g\partial \theta}\right |_{(e,0,1)}&0\\
\left .\frac {\partial ^2\Psi}{\partial h_1\partial \theta}\right
|_{(e,0,1)}&\cdots&
\left .\frac {\partial ^2\Psi}{\partial h_g\partial \theta}\right |_{(e,0,1)}&i&1\\
0&\ldots&0&1&0\end{array}\right ].$$ Therefore,
\begin{eqnarray}\det \left (H(\Psi)\right )=
\det \left (\left [\begin{array}{ccc} \left .\frac {\partial
^2\Psi}{\partial h_1^2}\right |_{(e,0,1)}&\cdots&
\left .\frac {\partial ^2\Psi}{\partial h_1\partial h_g}\right |_{(e,0,1)}\\
\vdots&\ddots&\vdots\\
\left .\frac {\partial ^2\Psi}{\partial h_g\partial h_1}\right
|_{(e,0,1)}& \cdots& \left .\frac {\partial ^2\Psi}{\partial
h_g^2}\right |_{(e,0,1)}\end{array}
\right ]\right )\nonumber \\
=\det \left (\left [
\begin{array}{ccc}\left .
\frac {\partial ^2\tilde \psi}{\partial h_1^2} \right
|_e&\cdots&\left .\frac {\partial ^2\tilde \psi}{\partial
h_1\partial h_g}
\right |_e\\
\vdots&\ddots&\vdots\\
\left .\frac {\partial ^2\tilde \psi}{\partial h_g\partial
h_1}\right |_e&\cdots& \left .\frac {\partial ^2\tilde
\psi}{\partial h_g^2}\right |_e
\end{array}
\right ]\right ),\label{eqn:hpsi}
\end{eqnarray}
where $\tilde \psi :U\rightarrow \mathbb{C}$ is the function
$g\mapsto \psi \left (\tilde \mu _{g^{-1}}(x),x\right )$. That
$(e,0,1)$ is a nondegenerate critical point then follows from the
following
\begin{lem} The Hessian of $\tilde \psi$ at $e$,
$H(\tilde \psi)_e$, is nonsingular. \label{lem:hessiane}
\end{lem}

\noindent {\it Proof.} We first produce an appropriate set of
local holomorphic coordinates on $M$ in the neighbourhood on
$p=\pi (x)$ and on $L^*$ in the neighbourhood of $x$, in terms of
which the action of $G$ will be a translation. Let $M_0=:\Phi
^{-1}(0)/G$ be the symplectic reduction of $M$. Then $M_0$ is an
$(n-g)$-dimensional complex manifold and has an induced K\"{a}hler
structure. Let $\alpha :\Phi ^{-1}(0)\rightarrow M_0$ be the
projection (a principal $G$-bundle) and set $\overline p=\alpha
(p)$. We claim:

\begin{lem}
There exist an open neighbourhood $T$ of $\overline p$ in $M_0$
and a section $\sigma :T\rightarrow \Phi ^{-1}(0)$ of $\alpha$,
such that $\sigma (\overline p)=p$ and $\sigma$ is holomorphic as
a map $T\rightarrow M$.\label{lem:section0}\end{lem}

\noindent \textit{Proof}. For $q\in M$, let $F_q\subseteq
T_qM\otimes \mathbb{C}$ be the $+i$-eigenspace of the complex
structure $J_q\in \mathrm{End}(T_qM)$.

Next, if $p\in \Phi ^{-1}(0)$, recall that $T_p\left (\Phi
^{-1}(0)\right )\subseteq T_pM$ is a coisotropic subspace, with
symplectic complement $T_p\left (\Phi ^{-1}(0)\right )^\perp
=T_p\left (G\cdot p\right )$ \cite{gs-gq}.

By assumption, $\Phi$ is submersive at $p\in \Phi ^{-1}(0)$. Let
$B$ be an open neighbourhood of $p$ on which $\Phi$ is a
submersion and the action of $G$ is locally free. Then $C=:\Phi
(B)$ is an open neighbourhood of $0\in \frak{g}^*$. For $c\in C$,
let $W_c=\Phi ^{-1}(c)\cap B$. Then $W_c$ is a
$(2n-g)$-dimensional real submanifold of $B$. For $q\in B$ let us
set $R_{q}=:T_{q}\left (W_{\Phi (q)}\right )^\perp$, the
symplectic complement of $T_{q}\left (W_{\Phi (q)}\right )$.

By Lemma 3.6 of \cite{gs-gq}, we have $F_p\cap \left (R_p\otimes
\mathbb{C}\right )=0$; therefore, perhaps after restricting $B$ to
a smaller open neighbourhood of $p$, we have $F_q\cap \left
(R_q\otimes \mathbb{C}\right )=0$ for all $q\in B$.

It follows that $F'_q=:F_q\cap \left (T_{q}\left (W_{\Phi
(q)}\right )\otimes \mathbb{C}\right )$ is a complex distribution
on $B$, of complex rank $n-g$. Being the intersection of two
integrable distributions, it is itself integrable.

Thus we may apply the complex Frobenius integrability theorem:
there are local holomorphic coordinates $r_i$ on $B\subseteq M$
centered at $p$ such that $F'=\{\partial /\partial
r_j:j=1,\ldots,n-g\}$. Lemma \ref{lem:section0} follows from this.

\bigskip

Let us now return to the proof of Lemma \ref{lem:hessiane}. Let
$z_1,\ldots,z_{n-g}$ be local holomorphic coordinates on $M_0$
centered at $\overline p$, and defined on an open neighbourhood
$T\subseteq M_0$ of $\overline p$. We may assume that the $z_i$'s
are induced by a biholomorphic diffeomorphism $\chi
:B_{n-g}(0,1)\rightarrow T$, where $B_m(w,r)$ denotes the ball
centered at the origin and of radius $r>0$ in $ \mathbb{C}^m$.

Let $\frak{g}_c=:\frak{g}\otimes \mathbb{C}$ be the Lie algebra of
the complexification $\tilde G$ of $G$. Let $S\subseteq
\frak{g}_c$ be an open neighbourhood of $0$ on which the
exponential map $\exp _{\tilde G}:\frak{g}_c\rightarrow \tilde G$
restricts to a biholomorphic diffeomorphism $S\rightarrow S'=:
\exp _{\tilde G}(S)$. Upon choosing an appropriate basis of
$\frak{g}$, which identifies $ \frak{g}_c$ with $ \mathbb{C}^g$,
we may assume that $S$ gets identified with $B_g(0,1)$. We now
define a holomorphic chart on $M$ on an open neighbourhood $D$ of
$p$ by setting
$$\gamma :B_{n-g}(0,1)\times B_g(0,1)\subseteq
\mathbb{C}^n\longrightarrow M,\,\,\,\,\,\,(w,z)\mapsto \exp
_{\tilde G}(w)\cdot \sigma \left (\chi (z)\right ).$$

Let us write $w=a+ib$ and $z=c+id$, for $w\in \mathbb{C}^{g}$ and
$z\in \mathbb{C}^{n-g}$, with $a,b\in \mathbb{R}^g$ and $c,d\in
\mathbb{R}^{n-g}$, and view $a,b,c,d$ as real local coordinates on
$M$ in the neighbourhood of $p$ . If $\xi \in \frak{g}$ has
coordinates $a=(a_1,\ldots,a_g)^t\in \mathbb{R}^g$ in the chosen
basis of $\frak{g}$, we have
$$\gamma ^{-1}\left (\exp _{G}(t\, \xi)\cdot p\right )=(ta,0).$$

We now construct local holomorphic coordinates on $L^*$ centered
at $x$. By $G$-invariance, $L,\,h_L,\,\nabla _L$ descend to
objects $L_0,\, h_0,\, \nabla _0$ on $M_0$. More precisely, with
obvious notation, if $\iota :\Phi ^{-1}(0)\hookrightarrow M_0$ is
the inclusion, we have $\iota ^*(L,h_L,\nabla _L)=\alpha
^*(L_0,h_0,\nabla _0)$. We choose a local holomorphic frame
$e_{L_0}$ for $L_0$ at $\overline p$, which perhaps after
restriction we may assume defined on the open neighbourhood $T$ of
Lemma \ref{lem:section0}. We shall view this as a holomorphic
section of $L$ defined on the complex submanifold $\alpha (T)$ of
$M$, and may assume with this identification that $
e_{L_0}(\overline p)=x$. Using the holomorphic action of $\tilde
G$ on $L$, we may then extend $e_{L_0}$ to a section $e_L$ of $L$
over the open neighbourhood $D$ of $p$ in $M$ by setting
$$e_L\left (\exp _{\tilde G}(w)\cdot \sigma (\overline q)\right )=
\exp _{\tilde G}(w)\cdot e_{L_0}\big (\sigma (\overline q)\big
)\,\,\,\,\,\,(w\in S,\,\overline q\in T).$$ Let $e_L^*$ be the
dual frame. The choice of $e_L^*$ then induces a holomorphic chart
$$\tilde \gamma : B_{n-g}(0,1)\times B_g(0,1)\times
\mathbb{C}\longrightarrow L^*_D=:\left .L^*\right |_D.$$ We have,
with $\xi \in \frak{g}$ and $a\in \mathbb{R}^{g}$ as above,
$$\tilde \gamma ^{-1}\left (\exp _{\tilde G}(t\xi)\cdot x\right )=
(ta,0,1).$$

Let now $\beta =:\left \| e_L^*\right \|^2$, so that $\beta
^{-1/2}e_L^*$ is a unitary frame of $L^*$ over $D$. We have an
induced trivialization $$\hat \gamma :B_{n-g}(0,1)\times
B_g(0,1)\times S^1\longrightarrow X^*_D=:\left .X^*\right |_D,$$
and since the action of $G$ preserves the metric we still have,
with $\xi \in \mathfrak{g}$ and $a\in \mathbb{R}^{g}$ as above,
\begin{eqnarray}\hat \gamma ^{-1}\left (\exp _{\tilde G}(t\xi)\cdot x\right
)= (ta,0,1).\label{eqn:loccoord}\end{eqnarray}

We now estimate $\psi \left (\tilde \mu _{\exp
_G(-t\,\xi)}(x),x\right )$ to second order using the Taylor
expansion of $\psi$ at $(x,x)$, (\ref{eqn:loccoord}) and
(\ref{eqn:psi}). Since $\xi _X(x)=\sum _{t=1}^ga_t'\left
.\frac{\partial}{\partial a_i}\right |_x$ is horizontal in view of
(\ref{eqn:hor}) and because $x\in \Phi ^{-1}(0)$, it is
annihilated by the connection form $\alpha _x$. Thus
$d_{(x,x)}\psi \left ((a,0)\right )=0$, and we have
\begin{eqnarray*}\psi \left (\tilde \mu _{\exp _G(-t\,\xi )}(x),x\right )=
\psi (x,x)+\frac 12 t^2a^t\, H(\tilde \psi)_e\,a+O(t^3)\\
= \psi (x,x)+\frac 12 t^2\,H(\rho)_X^{(2,0)}\big (\xi _X(x)\big
)+O(t^3), \,\,\,\,\,\,(t\in \mathbb{R}).
\end{eqnarray*}
where
\begin{eqnarray*}H(\rho)_x^{(2,0)}=\frac 12
\sum _{i,j}\frac {\partial ^2\rho}{\partial z_i\partial z_j}
(x)\,dz_i\,dz_j.
\end{eqnarray*}
Thus Lemma \ref{lem:hessiane} may be rephrased as:
\begin{lem}
The quadratic form $\alpha _x\in {\rm Sym}^2 \left
(\frak{g}^*\otimes \mathbb{C}\right )$ given by
\begin{eqnarray}\alpha _x:\xi \in \frak{g}\mapsto
H(\rho)_x^{(2,0)}\big (\xi _X(x)\big )\label{eqn:20}\end{eqnarray}
is nondegenerate (clearly $\xi _{L^*}(x)=\xi _X(x)$ when $x\in
X$). \label{lem:alfa} \label{lem:nondeg}\end{lem}

\noindent {\it Proof of Lemma \ref{lem:alfa}.} Let $H(\rho)_x$ be
the Hessian of $\rho$ at $x\in X\subset L^*$, in local holomorphic
coordinates induced by some holomorphic trivialization of $L$ in
the neighbourhood of $p=\pi (x)$; we may decompose it into types,
as
$$H(\rho)_x=H(\rho)_x^{(2,0)}+H(\rho)_x^{(1,1)}+H(\rho)_x^{(0,2)},$$
where
\begin{eqnarray*}
H(\rho)_x^{(1,1)}=\sum _{i,j} \frac {\partial ^2\rho}{\partial
z_i\partial \overline z_j} (x) \,dz_i\,d\overline z_j,\,\,\,
H(\rho)_x^{(0,2)}=\frac 12 \sum _{i,j}\frac {\partial
^2\rho}{\partial \overline z_i\partial \overline z_j}
 (x)\,d\overline z_i\,d\overline z_j.\end{eqnarray*}
Since $\rho$ is $G$-invariant we have, for any $\xi \in \frak{g}$:
\begin{eqnarray}0=H(\rho)_x\big (\xi _X(x)\big )=H(\rho)_x^{(2,0)}\big (\xi _X(x)\big )
+H(\rho)_x^{(1,1)}\big (\xi _X(x)\big )
\nonumber \\
+H(\rho)_x^{(0,2)}\big (\xi _X(x)\big )=H(\rho)_x^{(1,1)}\big (\xi
_X(x)\big )+ 2\,{\rm Re}\left (H(\rho)^{(2,0)}_x\big (\xi
_X(x)\big )\right ),\label{eqn:nondeg}
\end{eqnarray}
as $\rho$ is real-valued. By (\ref{eqn:hor}), if $x\in \Phi
^{-1}(0)$ then $\xi _X(x)$ is horizontal, that is, it lies in the
maximal complex subspace $H_x(X/M)$ of $T_x(X)\subset T_x(L^*)$.
By the ampleness of $L$, the Levi form $L_x(\rho)$ induces a
positive definite Hermitian form on $H_x(X/M)$. In turn, because
$G$ acts freely on $\Phi ^{-1}(0)$, the induced real quadratic
form on $\frak{g}$,
$$\xi \mapsto L_x(\rho)\big (\xi _X(x)\big )=:
H(\rho)_x^{(1,1)}\big (\xi _X(x)\big ),$$ is positive definite.
This proves Lemma \ref{lem:alfa}, because a complex symmetric
matrix with negative definite real part is nondegenerate.

\bigskip

Hence $(e,0,1)$ is a nondegenerate critical point of $\Psi$. We
still have to prove that, perhaps after restricting $U$ to a
smaller open neighbourhood of $e\in G$, it is the only critical
point of $\Psi$ on $U\times S^1\times (0,\infty)$.

Now $(0,1)$ is the only critical point of $\Psi (x,e,\theta,u)$ as
a function of $(\theta,u)$. Therefore, given any $\epsilon >0$ we
may replace $U$ with some possibly smaller open neighbourhood of
$e$, so that $(e,0,1)$ is the only critical point of $\Psi$ in
$U\times S^1\times (1-\epsilon,1+\epsilon)$. Suppose that there is
a sequence $(g_i,\theta_i,t_i)$ of critical points of $\Psi$ with
$g_i\rightarrow e$. Since $S^1$ is compact, the sequence $\theta
_i$ has some accumulation point $\theta _\infty$, and after
passing to a subsequence this implies $(d\Psi)_{(e,\theta
_\infty,t_i)}\rightarrow 0$. By (\ref{eqn:critpsi}), we must have
$\theta _\infty=0$ and  $t_i\rightarrow 1$, against the fact that
$|t_i-1|\ge \epsilon$.

This completes the proof of the Lemma.

\bigskip

We are now in a position to apply the complex stationary phase
Lemma to estimate the asymptotic behaviour of (\ref{eqn:komega2})
as $k\rightarrow +\infty$ \cite{hor}, \cite{ms}. More precisely,
since $\chi _\omega (e)=\dim (V_\omega)$,
\begin{eqnarray*}\label{eqn:asympt}
\Pi _{k,\omega}(x,x)'\sim k^{n-g/2}\, \dim (V_\omega)^2\,
\det \left (H(\Psi)(e,0,1)/2\pi i\right )^{-1/2}\, s_0(x,x)\\
+ O\left (k^{n-g/2-1}\right ).\end{eqnarray*} Theorem \ref{thm:1}
follows in view of (\ref{eqn:s0}), (\ref{eqn:hpsi}) and Lemma
\ref{lem:hessiane}.

\section{Proof of Theorem \ref{thm:3}.} We have
\begin{eqnarray}\label{eqn:integral}
\dim H^0(M,L^{\otimes k})_\omega \,=\,\int _X\, \Pi _{\omega
,k}(x,x)\, dx,\end{eqnarray} and we want to estimate
asymptotically the latter integral as $k\rightarrow +\infty$.

As we have seen, the equivariant Szeg\"{o} kernel $\Pi _{\omega
,k}$ is given by
\begin{eqnarray*}
\Pi _{\omega ,k}(x,y)\, =\frac{1}{2\pi}\, \int _0^{2\pi}\, \int _G
\,e^{-ik\theta}\,\chi _\omega \left (g^{-1 }\right)\, \Pi \left
(e^{i\theta}\mu _{g^{-1}}(x),y\right )\, d\theta\,
dg\,\,\,\,\,\,(x,y\in X).
\end{eqnarray*}
Here $\mu$ is the action of $G$ on $X$, and $\rho$ is induced by
$\mu$ by the formula $\rho _g (f)(x)=f\big (\mu _{g^{-1}}(x)\big
)$ ($x\in X$, $g\in G$).

By compactness and the results of \cite{bs}, we may find
$\varepsilon _0>0$ with the following property: For any $p\in M$
and $\varepsilon >0$, let $B(p,\varepsilon )$ be the open ball in
$M$ centered at $p$ and radius $\varepsilon $, in the geodesic
distance on $M$. Then on the inverse image $\big (\pi \times \pi
\big ) ^{-1}\left (B(p,2\varepsilon _0)\times B(p,2\varepsilon
_0)\right )$ the Szeg\"{o} kernel is microlocally equivalent to a
Fourier integral operator of the type (\ref{eqn:fourierintegral}).

We shall denote by $\mathrm{dist}_M$ the geodesic distance
function on $M$, as well as its pull-back to $X$, and by
$\mathrm{dist}_X$ the geodesic distance function on $X$. Here the
Riemannian metric on $X$ is defined in the natural manner in terms
of the metric on $M$ and the connection, so that the projection
$\pi :X\rightarrow M$ is a Riemannian submersion;
$\mathrm{dist}_X$ is clearly $S^1$-invariant. Let us set:
\begin{eqnarray}
\mathcal{V}=\left \{(x,g)\in X\times G\, :\,\mathrm{dist}_M\left
(x,\mu _{g^{-1}} (x)\right )<\epsilon _0^2\right
\}.\label{eqn:calv}
\end{eqnarray}
$\mathcal{V}$ is an $S^1$-invariant open neighbourhood of $X\times
\{e\}$ in $X\times G$. We now decompose the integral
(\ref{eqn:integral}) as follows:
\begin{eqnarray}
\lefteqn{\int _X\, \Pi _{\omega ,k}(x,x)\, dx=}\nonumber \\
& & =\, \frac{1}{2\pi}\int _X\int _G\int _0^{2\pi}\,
\,e^{-ik\theta}\,\chi _\omega \left (g^{-1 }\right)\, \Pi \left
(e^{i\theta}\mu _{g^{-1}}(x),x\right )\, d\theta\, dg\, dx\nonumber \\
& & =\, \int \!\int _{\cal V}\, \chi _\omega \left (g^{-1
}\right)\, \Pi _k\left (\mu _{g^{-1}}(x),x\right )\, dx\, dg\,
\nonumber \\
& & +\, \int \!\int _{X\times G\setminus {\cal V}}\, \,\chi
_\omega \left (g^{-1 }\right)\, \Pi _k\left (\mu
_{g^{-1}}(x),x\right
)\,dx\, dg    \nonumber \\
& & =H^{(\omega ,k)}_1+H_2^{(\omega,k)}.\label{eqn:decompose}
\end{eqnarray}
$H^{(\omega ,k)}_1$ and $H_2^{(\omega,k)}$ are defined by the
latter equality. $H_2^{(\omega,k)}$ is $O(k^{-N})$ for every
$N=1,2,\ldots$, since the Szeg\"{o} kernel is smoothing away from
the diagonal. We thus need to estimate $H^{(\omega ,k)}_1$
asymptotically as $k\rightarrow +\infty$.

To this end, we shall first of all construct an appropriate finite
open cover of $M$. Choose points $p_j\in \Phi ^{-1}(0)$ ($1\le
j\le r_0$ for some fixed $r_0\ge 1$) such that
$$ \Phi ^{-1}(0)_{\varepsilon _0/2} \subseteq
\bigcup _{j=1}^{r_0} B(p_j,\varepsilon _0) ;$$ here $\Phi
^{-1}(0)_{\varepsilon _0/2}$ is the $\varepsilon
_0/2$-neighbourhood of $\Phi ^{-1}(0)$ in $M$ in the geodesic
distance.

Let us further choose points $p_j\in M\setminus \Phi
^{-1}(0)_{\varepsilon _0/2}$ ($r_0+1\le j\le r_1$ for some integer
$r_1>r_0$), such that
$$M\setminus \Phi ^{-1}(0)_{\varepsilon
_0/2}\subseteq \bigcup _{j=r_0+1}^{r_1} B(p_j,\varepsilon _0/3).$$
Clearly, $r_0$ and $r_1$ depend on $\epsilon _0$.

Let us then set $V_j=B(p_j,\varepsilon _0)$ if $1\le j\le r_0$,
$V_j=B(p_j,\varepsilon _0/3)$ if $r_0+1\le j\le r_1$. This is a
finite open cover of $M$. Let $1 =\sum _j\phi _j$ be a partition
of unity subordinate to the cover $\{V_j\}$. We shall also write
$\phi _j$ for $\phi _j\circ \pi$. Thus, $1=\sum _j\phi _j$ will be
implicitly seen as a partition of unity on $X$ subordinate to the
open cover $\tilde V_j$, where $\tilde V_j=\pi ^{-1}(V_j)$.

$H_1^{(\omega ,k)}$ may be decomposed as
\begin{eqnarray}\label{eqn:integraldecomp}
H_1^{(\omega ,k)}\,=\sum _j\, \int \!\int _{{\cal V}_j}\, \chi
_\omega \left (g^{-1 }\right)\, \phi _j(x)\,\Pi _k\left (\mu
_{g^{-1}}(x),x\right )\, dx\, dg =\sum _j
A^{(\omega,k)}_j,\end{eqnarray} where $\mathcal{V}_j=\left
\{(x,g)\in \mathcal{V}:x\in \tilde V_j\right \}$ and
$A^{(\omega,k)}_j$ is defined by the latter equality. We shall
next estimate each $A^{(\omega,k)}_j$ separately.

By construction, if $(x,g)\in \mathcal{V}$ then the neighbourhood
of $\left (\mu _{g^{-1}}(x),x\right )$ we may represent $\Pi$ as a
Fourier integral of the type (\ref{eqn:fourierintegral}).

More precisely, by our choice of $\varepsilon _0$, for every
$j=1,\ldots,r_1$ we may choose a preferred holomorphic section of
$L$ at $p_j$ and adapted local holomorphic coordinates at $p_j$,
in the sense of \cite{sz}, both defined on $B(p_j,\varepsilon
_0)$. We then have induced holomorphic coordinates on $L$. Let
$a_j=a_j(z)$ be the square norm in this preferred frame, and let
$a_j(z,w)$ be its extension to $M\times M$, almost analytic in
$z$, and almost antianalytic in $w$. This determines a phase
function $\hat \psi _j(z,\lambda,w,\mu)=i\left (1-a_j(z,w)\lambda
\overline \mu \right )$, restricting to $\psi _j=i\left
(1-\frac{a_j(z,w)}{\sqrt{a_j(z)}\sqrt{a_j(w)}}\lambda \overline
\mu \right )$ on $X\times X$ (with $\lambda, \mu \in S^1$). $\psi
_j$ locally parametrizes the almost analytic Lagrangian relation
associated to the Szeg\"{o} kernel. In particular,
$$d_{(x,x)}\psi =\left (\partial _z\psi (x,x),
\overline \partial _w\psi (x,x)\right )=(\alpha _x,-\alpha _x),$$
for every $x\in \tilde V_j$. A straighforward computation then
gives
$$d_{(e^{i\theta}x,x)}\psi =(e^{i\theta}\alpha
_{e^{i\theta}x},-e^{i\theta}\alpha _x)\,\,\,\,\,\, (x\in \tilde
V_j,\,e^{i\theta}\in S^1).$$ Let $s$ be the classical symbol
appearing in the local representation of $\Pi$ as a Fourier
integral operator. Setting $s_j(y,x,t)=\phi _j(x)\, s_j(y,x,t)$,
we obtain
\begin{eqnarray}\label{eqn:H3j}
A_j^{(\omega,k)}= \frac{1}{2\pi}\, \int \! \int
_{\mathcal{V}_j}\,\int _0^{2\pi}\! \int _0^{+\infty}
\,e^{-ik\theta}\,\chi _\omega \left (g^{-1 }\right)e^{it\psi
_j\left (e^{i\theta}\mu _{g^{-1}}(x),x\right )}
 \nonumber \\
 \times s_j\left (e^{i\theta}\mu _{g^{-1}}(x),x,t \right )\, dx\,
d\theta \,dg\, dt.
\end{eqnarray} With the change of variables $t=ku$, this is
\begin{eqnarray}\label{eqn:H3jbis}
A_j^{(\omega,k)} =    \frac{k}{2\pi}\, \int \! \int
_{\mathcal{V}_j}\!\int _0^{2\pi}\! \int _0^{+\infty} \,e^{ik\left
( u\psi _j(e^{i\theta}\mu _{g^{-1}}(x),x)- \theta \right )} \,\chi
_\omega \left (g^{-1 }\right)\, &
\nonumber \\
\times \, s_j\left (e^{i\theta}\mu _{g^{-1}}(x) ,x,ku\right ) \,
dx\, d\theta \,dg\,
 du &
\nonumber \\
=\frac{k}{2\pi}\, \int \!\int _{\mathcal{V}_j}\int _0^{2\pi}\!\int
_0^{+\infty} \,e^{ik \Psi _j(x,g,\theta,u) } \,\chi _\omega \left
(g^{-1 }\right) \,
 &\nonumber \\
\times s_j\left ( e^{i\theta}\mu _{g^{-1}}(x) ,x,ku\right ) \,
dx\, d\theta \, du\, dg,&
\end{eqnarray}
where $\Psi _j:\mathcal{V}_{j}\times S^1\times
\mathbb{R}_+\rightarrow \mathbb{C}$ is defined by
\begin{eqnarray}\label{eqn:phase}
\Psi _j(x,g,\theta,u) \, = \, u \psi _j(e^{i\theta}\mu
_{g^{-1}}(x),x) - \theta .
\end{eqnarray}
\noindent To estimate $A_j^{(\omega,k)}$, we shall now use the
complex stationary phase Lemma. For $x\in \tilde V_j$, let $\Psi
_x:U_x\times S^1\times (0,+\infty)\rightarrow \mathbb{C}$ be the
partial function $\Psi _{jx}(g,\theta ,u)=:\Psi _j(x,g,\theta,u)$;
here $U_x\subseteq G$ is a suitable neigbourhood of $e\in G$.

\begin{lem}
For all sufficiently small $\varepsilon _0$, there exists $\gamma
_0>0$ (independent of $\epsilon _0$) such that $\left \|d\Psi
_{jx}\right \| \ge \gamma _0\, \epsilon _0$ at any
$(x,g,\theta,u)\in \mathcal{V}_j\times S^1\times \mathbb{R}_+$
with $ r_0+1\le j\le r_1 $. \label{lem:psix}\end{lem}

\noindent \textit{Proof.} To begin with, we make the following
remark. Given any $a,b$ with $0<a<1<b$, there exists $d>0$ with
the following property: For all sufficiently small $\epsilon _0>0$
and $(x,g,\theta,u)$ with $(x,g)\in \mathcal{V}$, we have $\left
\|d\Psi _x\right \|>d$ if $u\not \in (a,b)$.

If not, there would exist a sequence $(x_i,g_i,\theta _i,u_i)$
with
$$\mathrm{dist}_M\left (\mu _{g_i^{-1}}(x_i),x_i\right
)\rightarrow 0,$$ $u_i\not \in (a,b)$, and $d_{(x_i,g_i,\theta
_i,u_i)}\Psi _{jx}\rightarrow 0$. Given the compactness of $X$ and
$G$, we may pass to a subsequence and assume $x_i\rightarrow
x_\infty$, $g_i\rightarrow g_\infty$. Then $\mu _{g_\infty
^{-1}}(x_\infty )=e^{i\vartheta _\infty}x_\infty$, for some
$\vartheta _\infty \in [0,2\pi)$, and $d_{(x_\infty,e^{i\vartheta
_\infty}x_\infty,\theta _i,u_i)}\Psi _{jx}\rightarrow 0$. But by
the analysis in \cite{z} of the critical points of $\psi$ this
implies $\theta _i\rightarrow -\vartheta _\infty$, $u_i\rightarrow
1$, absurd. Thus, we may assume $u\in [a,b]$.

We shall then fix some $a\in (0,1)$ and suppose that $\epsilon _0$
is sufficiently small for the conclusion of Lemma \ref{lem:psix}
to hold.

If $(x,g)\in \mathcal{V}_j$, let $\vartheta =\vartheta (x,g)\in
[0,2\pi)$ be uniquely determined by the condition that
$$\mathrm{dist}_X\left (e^{i\vartheta}\mu _{g^{-1}}\left (x\right
),x\right )=\mathrm{dist}_M\left (\mu _{g^{-1}}\left (x\right
),x\right )<\epsilon _0^2.$$ Given that
$d_{(e^{i\theta}x,\,x)}\psi _j=(e^{i\theta}\alpha
_{e^{i\theta}x},-e^{i\theta}\alpha _x)$ for all $x\in \tilde V_j$,
there exists $C_1>0$ such that
\begin{equation}(x,g)\in \mathcal{V}_j\, \Longrightarrow \, \left
\| d_{(\mu _{g^{-1}}(e^{i\theta}x),\,x)}\psi _j-(e^{i(\theta
+\vartheta)}\alpha _{\mu
_{g^{-1}}(e^{i(\theta+\vartheta)}x)},-e^{i(\theta+\vartheta)}\alpha
_x)\right \|<C_1\, \epsilon
_0^2.\label{eqn:estimateondpsi}\end{equation}

Suppose now that $r_0+1\le j\le r_1$ and $p\in V_j$. Then
$${\rm dist}\left (p, \Phi ^{-1}(0)\right ) \ge {\rm dist}\left (p_j,
\Phi ^{-1}(0)\right ) - d(p,p_j)\ge \varepsilon _0/6.$$ Suppose
$x\in \tilde V_j$ and $(x,g)\in \mathcal{V}$, and set $p=\pi
(x)\in V_j$. Then $$\mathrm{dist}_M(\mu _{g^{-1}}(p),p)=
\mathrm{dist}_M(\mu _{g^{-1}}(x),x)<\varepsilon _0^2.$$ Therefore,
$$\mathrm{dist}_M\left (\mu _{g^{-1}}(p),\Phi ^{-1}(0)\right ) \ge
\varepsilon _0/6-\varepsilon _0^2\ge \varepsilon _0/7,$$ if
$\varepsilon _0<1/42$. Thus, given that $0\in \frak{g}^*$ is a
regular value of the moment map, there exists $C_2>0$ such that
\begin{equation} \label{eqn:moment} r_0+1\le j\le r_1,\, (x,g)\in
\mathcal{V}_j\, \Longrightarrow \left \| \Phi \left (\mu
_{g^{-1}}(p)\right )\right \|\ge C_2\,\varepsilon
_0,\end{equation} in a given fixed norm on $\frak{g}^*$. In other
words, if $r_0+1\le j\le r_1$ and $(x,g)\in \mathcal{V}_j$ then
there exists $\xi \in \frak{g}^*$ of unit norm such that $\Phi
_\xi \left (\mu _{g^{-1}}(p)\right )\ge C_2\,\varepsilon _0$,
where $p=\pi (x)$ and $\Phi _\xi =\big <\Phi ,\xi \big >$. This is
equivalent to the condition \begin{equation} \alpha _{\mu
_{g^{-1}}(x)}\left (\xi _X\left ( \mu _{g^{-1}}(x)\right )\right
)\ge C_2\, \varepsilon _0,\label{eqn:estimateonalfa}\end{equation}
and in view of the above this completes the proof of Lemma
\ref{lem:psix}.

\bigskip

We now apply the complex stationary phase Lemma to estimate the
sum $\sum _{j=r_0+1}^{r_\epsilon}H_{3,j}^{(\omega,k)}$. Recalling
the definition of $H_{3,j}^{(\omega,k)}$, we first integrate over
$\theta,\,u,\,g$ and then over $x$. The first integral, given
Lemma \ref{lem:psix}, is $\le C_Nk^{-N}$ for every $N=1,2,\ldots$.
We thus obtain
$$\left |\sum _{j=r_0+1}^{r_\epsilon}H_{3,j}^{(\omega,k)}\right
|\le C_Nk^{-N}\int _X\left (\sum _{j=r_0+1}^{r_\epsilon}\phi
_j\right )\,dx\le C_N'k^{-N}.$$

In order to estimate $\sum _{j=1}^{r_0}H_{3,j}^{(\omega,k)}$, we
shall use in the neighbourhood of any $x\in \pi^{-1}\left (\Phi
^{-1}(0)\right )$ the local coordinates $(w,z,\lambda
=e^{i\theta})$ discussed in the proof of Lemma \ref{lem:hessiane}.
Thus $w\in \mathbb{C}^g$, $z\in \mathbb{C}^{n-g}$ and $\lambda \in
S^1$. We shall write $w=a+ib$ and $z=c+id$, where $a,b\in
\mathbb{R}^g$ and $c,d\in \mathbb{R}^{n-g}$. Then $\pi ^{-1}\left
(\Phi ^{-1}(0)\right )$ is locally defined near $x$ by the
equation $\{b=0\}$. In the neighbourhood of $(x,e,1,1)\in X\times
G\times S^1\times \mathbb{R}$ we then have local coordinates
$(a,b,c,d,\lambda =e^{i\theta}, a',\lambda '=e^{i\theta '},u)$. In
the following we shall not distinguish between $\Psi$ and its
expression in local coordinates.

Loosely speaking, $\pi ^{-1}\left (\Phi ^{-1}(0)\right )\times
\{(e,1,1)\}$ is a nondegenerate critical manifold for the (locally
defined) phase $\Psi$.

\begin{lem}
Every $x\in \pi ^{-1}\left (\Phi ^{-1}(0)\right )\times
\{(e,1,1)\}$ is a critical point of $\Psi$. The Hessian of $\Psi$
at any $x\in \pi ^{-1}\left (\Phi ^{-1}(0)\right )\times
\{(e,1,1)\}$ has rank $2g+2$, and is nondegenerate in the
variables $(b,a',\theta ',u)$.\label{lem:critical1}\end{lem}

\noindent \textit{Proof.} By the previous discussion, $\pi
^{-1}\left (\Phi ^{-1}(0)\right )\times \{(e,1,1)\}$ is a critical
manifold for $\Psi$. Thus the Hessian of $\Psi$ at any $x\in \pi
^{-1}\left (\Phi ^{-1}(0)\right )\times \{(e,1,1)\}$, $H(\psi)_x$,
is a well-defined quadratic form on the tangent space to $X\times
G\times S^1\times \mathbb{R}$ at $x$. Suppose without loss of
generality that the local coordinates $(a,b,c,d,\lambda
=e^{i\theta })$ are centered at $x$. Clearly, every second
derivative of $\Psi$ involving one of the variables $a,c,d,\theta$
(which give local coordinates on $\pi ^{-1}\left (\Phi
^{-1}(0)\right )\times \{(e,1,1)\}$) vanishes at $x$
(incidentally, given the expression in local coordinates of the
phase $\psi$ discussed in \cite{z}, $\Psi$ does not depend on
$\lambda$). We want to show that $H(\psi)_x$ is nondegenerate on
the subspace
$${\rm span}\left \{\left .\frac{\partial}{\partial
 b_i}\right |_x, \left .\frac{\partial}{\partial a'_i}\right
|_x, \left .\frac{\partial}{\partial \theta'}\right |_x, \left
.\frac{\partial}{\partial u}\right |_x\right \}.$$ Let us then
consider then the $(2g+2)\times (2g+2)$ matrix given by Hessian of
$\Psi$ at $x$ in the variables $(b,a',\theta,u)$, $\hat
H(\psi)_x$. To begin with, let us remark that
\begin{claim}
$$\left . \frac{\partial ^2\Psi}{\partial u\,\partial a'}\right |_x=0,\,
\left .\frac{\partial ^2\Psi}{\partial u\,\partial b}\right
|_x=0.$$
\end{claim}

\noindent {\it Proof.} By definition, $\frac{\partial
\Psi}{\partial u}=\psi$. The first vanishing then holds because
 $d_{(x,x)}\psi
=(\alpha _x,-\alpha _x)$, where $\alpha$ is the connection 1-form,
and as we have seen the action of $G$ at $x\in \pi ^{-1}\left (
\Phi ^{-1}(0)\right )$ is horizontal to first order. The second
vanishing holds because $\psi$ is constant along the diagonal.

As recalled in the proof of Theorem \ref{thm:1}, the Hessian of
$\Psi$ in the two variables $\theta ',\,u$ at $(0,1)$ is known
from \cite{z}. Given the Claim, we obtain for the Hessian in
$b,a',\theta ',u$:
\begin{eqnarray*}
\det \hat H(\Psi)_x=\det \left (\left [\begin{array}{ccc} \left .
\frac{\partial ^2\Psi }{\partial  b^2}\right |_x&\left .
\frac{\partial ^2\Psi}{\partial a'\, \partial b}\right |_x
\\ \left .\frac{\partial ^2\Psi}{\partial a'\, \partial b}\right |_x&
\left .\frac{\partial ^2\Psi}{\partial a^{\prime 2}}\right |_x
\end{array}\right ]\right )=\det \left (\left [\begin{array}{ccc} \left .
\frac{\partial ^2\psi }{\partial  b^2}\right |_x&\left .
\frac{\partial ^2\psi}{\partial a'\, \partial b}\right |_x
\\ \left .\frac{\partial ^2\psi}{\partial a'\, \partial b}\right |_x&
\left .\frac{\partial ^2\psi}{\partial a^{\prime 2}}\right |_x
\end{array}\right ]\right ).
\end{eqnarray*}
\begin{claim} $$\left .\frac{\partial ^2\Psi }{\partial
b^2}\right |_x=0.$$
\end{claim}

\noindent \textit{Proof.\textit{}} If we fix $g=e,\theta '=0$,
$u=1$, $\Psi =\psi$ is constant as a function of $x$.

Thus,
\begin{eqnarray*}
\det \hat H(\Psi)_x=\det \left (\left [\begin{array}{ccc} 0&\left
.\frac{\partial ^2\psi}{\partial a' \partial b}\right |_x\\
\left .\frac{\partial ^2\psi}{\partial a'
\partial b}\right |_x& \left .
\frac{\partial ^2\psi }{\partial  a^{\prime 2}}\right
|_x\end{array}\right ]\right ),
\end{eqnarray*}
and we are reduced to proving that the $g\times g$ symmetric
matrix $ \left .\frac{\partial ^2\Psi}{\partial a'\, \partial
b}\right |_x=\left .\frac{\partial ^2\psi}{\partial a'\, \partial
b}\right |_x$ is nonsingular. We shall do this using the second
order expansion of $\Psi =\psi$ as a function of $a'$ and $b$
keeping fixed $a=0$, $\theta =\theta '=0$, $u=1$. In so doing, we
shall use the Taylor expansion of $\psi$ (as a function on $L^*$)
at $(x,x)$; this is determined as we have seen by the Taylor
expansion of the square norm function $\varrho$.

In the notation of the proof of Lemma \ref{lem:hessiane}, let
$\gamma (a+ib,c+id)$, $\tilde \gamma (a+ib,c+id,z)$ and $\hat
\gamma (a+ib,c+id,\lambda)$ be, respectively, the local charts of
$M$, $L^*$ and $X\subseteq L^*$ near $p=\pi (x)$ and $x$; here
$a,\, b\in \mathbb{R}^g$, $c,\, d\in \mathbb{R}^{n-g}$, $z\in
\mathbb{C}$, $\lambda \in S^1$. Thus,
$$\gamma (a+ib,0)=\exp (a+ib)\cdot p,$$
$$\tilde \gamma (a+ib,0,z)=\exp (a+ib)\cdot (p,z\, e_L^*(p)),$$
$$\hat
\gamma \left (a+ib,c+id,\lambda \right )=\tilde \gamma \left
(a+ib,c+id,\lambda /\sqrt{\beta \left (\gamma (a+ib,c+id)\right
)}\right ),$$ where $\beta =\left \|e_L^*\right \|^2$. To simplify
notation, we identify $a+ib\in \mathbb{C}^g$ with the
corresponding vector in $\frak{g}_c$, given the implicit choice of
a fixed basis of $\frak{g}$. We also write $g=\exp (a+ib)$ for its
image in $\tilde G$ under the exponential map. With some further
abuse of notation, in the following lines we shall omit the
variables $c,\, d$: for example, $\hat \gamma (i\,t\,b,\lambda)$
will then really mean $\hat \gamma (i\,t\, b,0,\lambda)$, where
$\lambda =e^{i\theta}$.

We shall also let $p=\pi (x)$, so that $\gamma (a+i\,b)=\exp
(a+i\,b)\cdot p$, and $x=(p,\tilde e_L^*(p))$. By construction of
$\tilde e_L$, if $a'\in \mathbb{R}^g$ and $t\in \mathbb{R}$ is
sufficiently small, we have
$$\exp (t\, a')\cdot \left (q,\lambda \, \tilde e_L^*(q)\right )=
\left (\exp (t\, a')\cdot q,\lambda \, \tilde e_L^*(\exp (t\,
a')\cdot q)\right ),$$ for every $q\in M$ near $p$ and $\lambda
\in \mathbb{C}$. In particular, if $q=\gamma (i\,t\,b)=\exp
(i\,t\,b)\cdot p$, we have
\begin{eqnarray}
e^{t\,a'}\cdot \tilde \gamma (i\,t\,b,\lambda)= e^{t\,a'}\cdot
\Big ( \gamma (i\,t\,b),\, \lambda \, \tilde e_L^*\big (\gamma
(i\,t\,b)\big )\Big )\nonumber \\= e^{t\, a'}\cdot \Big (
e^{i\,t\,b}\cdot p,\, \lambda \,e^{i\,t\,b}\cdot \tilde e_L^*\big
(e^{i\,t\,b}\cdot p \big )\Big ) \nonumber
\\
=e^{t\, a'}\cdot e^{i\,t\,b}\cdot \big (p,\, \lambda \, \tilde
e_L^*(p)\big )\nonumber \\
=\Big ( e^{t\,a'}e^{i\,t\,b}\cdot p,\, \lambda \, \tilde e_L^*\big
(e^{t\,a'}e^{i\,t\,b}\cdot p \big )\Big ).\end{eqnarray} Since
$\exp (t\, a')\exp (i\,t\,b)=\exp \Big (
t(a'+i\,b)+(i/2)\,t^2\,[a',b]+O(t^3)\Big )$, we may rewrite this
as
\begin{eqnarray*}
\lefteqn{ e^{t\,a'}\cdot \tilde \gamma (i\,t\,b,\lambda) =} \\
& & =
 \Big ( e ^{t\,(a'+i\,b)+(i/2)\,t^2\,[a',b]+O(t^3)}\cdot
p, \lambda \,\tilde e_L^*\big (e ^{t\,(a'+i\,b)+(i/2)
\,t^2\,[a',b]+O(t^3)}\cdot p \big )\Big )\nonumber
\\
& & =\tilde \gamma \left (
t\,(a'+i\,b)+(i/2)\,t^2\,[a',b]+O(t^3),\lambda \right
).\end{eqnarray*}

\noindent We have for $a',\, b\in \mathbb{R}^g$ and $t\in
\mathbb{R}$ sufficiently small:
\begin{eqnarray*}
\Psi \left (\hat \gamma ( i\,t\,b,1) ,e^{-t\,a'},0,1,\right )=
\Psi \left (\tilde \gamma \left ( i\,t\,b,\frac{1}{\beta (\gamma
(i\,t\,b))}\right ),e^{-t\,a'},0,1,\right )\nonumber \\
= \psi  \left ( e^{t\,a'}\cdot \left ( e^{i\,t\,b}\cdot
p,\frac{1}{\beta (\gamma (i\,t\,b) )} \right ), \left (
e^{i\,t\,b}\cdot p,\frac{1}{\beta (\gamma (i\,t\,b) )} \right )
\right ),\end{eqnarray*} Identifying $\psi$ with its expression in
local coordinates this is
\begin{eqnarray*}\psi \left ( \left
(t\,(a'+i\,b)+(i/2)\,t^2\,[a',b]+O(t^3),\, \frac{1}{\beta (\gamma
(i\,t\,b) )}\right ), \left (i\,t\,b,\frac{1}{\beta (\gamma
(i\,t\,b) )}\right )\right ).
\end{eqnarray*}
Now, as we have mentioned, the action of $G$ on $X\subseteq L^*$
is horizontal over $\Phi ^{-1}(x)$. Thus, given that $\pi (x)\in
\Phi ^{-1}(x)$, for every $\xi \in \frak{g}$ the associated
tangent vector $\xi _{L^*}(x)=\xi _{X}(x)\in H(X/M)_x$. The latter
is the maximal complex subspace of $T_xL^*$ contained in $T_xX$.
Therefore, we also have $J_x(\xi _{X}(x))=(i\xi )_X(x)\in
H(X/M)_x$. Now, in more intrinsic notation, suppose $a',\,b\in
\mathbb{R}^g$ correspond to $\xi ',\, \xi \in \frak{g}$,
respectively; then, in the local coordinates provided by the
holomorphic chart $\tilde \gamma$, $(a',0)$ stands  for $\xi
'_{L^*}(x)$, $(b,0)$ for $\xi _{L^*}(x)$, and $(ib,0)$ for $J_x\xi
_{L^*}(x)$ ($J_x$ is the complex structure on $T_xL^*$). Thus, on
the one hand in our local coordinates
$$\alpha _x\big (t\,(a'+i\,b)+(i/2)\,t^2\,[a',b],0\big
)=0,\,\,\alpha _x\big (i\,t\,b,0)=0.$$ On the other hand, the path
$\tilde \gamma (i\, t\,b)$ is tangent to $X$ at $t=0$, and
therefore $\beta (\gamma (i\,t\,b) )^{-1}=1+t^2\,s(t)$ for some
smooth function $s(t)$. Given that $d_{(x,x)}\psi =(\alpha
_x,-\alpha _x)$, where $\alpha$ is the connection 1-form, we
conclude from the latter expressions that
\begin{eqnarray*}
\Psi \left (\hat \gamma ( i\,t\,b,1) ,e^{-t\,a'},0,1,\right
)=\frac{1}{2}\,t^2\,H(\psi) _{(x,x)}\big ( a'+ib,ib\big
)+O(t^3).\end{eqnarray*} In view of (\ref{eqn:psi}),
\begin{eqnarray}
\lefteqn{H(\psi) _{(x,x)}\big ( a'+ib,ib\big
)=H^{(2,0)}(\varrho)_x(a'+ib,a'+ib)}\nonumber
\\& &+H^{(1,1)}(\varrho)_x(a'+ib,ib)
+H^{(0,2)}(\varrho)_x(ib,ib).\label{eqn:hessianrho}
\end{eqnarray}
Since $\psi (x,x)$ vanishes identically for $x\in X$,
\begin{eqnarray*}
\psi  \left ( \left ( \gamma (i\,t\,b) ,\frac{1}{\beta (\gamma
(i\,t\,b) )} \right ), \left ( \gamma (i\,t\,b) ,\frac{1}{\beta
(\gamma (i\,t\,b) )} \right ) \right )=0\end{eqnarray*} for every
$t$, and therefore in (\ref{eqn:hessianrho}) the terms containing
only $b$ add up to zero. Using the symmetry and the
$\mathbb{C}$-linearity of $H^{(2,0)}(\varrho)$ and the
sesquilinearity of $H^{(1,1)}(\varrho)$, we obtain
\begin{eqnarray}
\lefteqn{H(\psi) _{(x,x)}\big ( a'+ib,ib\big
)=} \nonumber \\
& & =H^{(2,0)}(\varrho)_x(a',a')+i\big (
2H^{(2,0)}(\varrho)_x-H^{(1,1)}(\varrho)_x\big
)(a',b).\label{eqn:hessianrho1}
\end{eqnarray}

Thus, the $g\times g$ complex matrix $-i\left . \frac{\partial
^2\Psi}{\partial a'\, \partial b}\right |_x$ represents the
bilinear pairing $\frak{g}\times \frak{g}\rightarrow \mathbb{C}$
given by
$$(\xi,\eta)\mapsto
2H^{(2,0)}(\varrho)_x\big (\xi _X(x),\eta _X(x)\big
)-H^{(1,1)}(\varrho)_x\big (\xi _X(x) ,\eta _X(x)\big ).$$ The
latter is in turn the restriction of a similarly defined
$\mathbb{R}$-bilinear pairing on $\frak{g}_c$. We have seen in the
proof of Theorem \ref{thm:1} that the former term is represented
by a complex symmetric matrix with negative definite real part,
while the second term is a negative definite Hermitian product on
$\frak{g}_c$. The proof of Lemma \ref{lem:critical1} is then
completed by the following:

\begin{lem} Let $b\ge 1$ be an integer.
Let $C$ and $H$ be complex $b\times b$ matrices. Suppose that $C$
is symmetric and $H$ is Hermitian. Suppose furthermore that $C$
has positive definite real part, and that $H$ is positive
semidefinite. Then $C+H$ is nonsingular.\label{lem:sh}
\end{lem}

\noindent \textit{Proof.} By the hypothesis, $C=A+iB$, $H=R+iS$,
where $A,B,R,S$ are real $b\times b$ matrices satisfying: $A^t=A$,
$B^t=B$, and $A$ is positive definite; $R^t=R$, $S^t=-S$ and
$$(V^t+iW^t)(R+iS)(V-iW)\ge 0$$
for all $V,W\in \mathbb{R}^b$, and equality holds if and only if
$V=W=0$. Suppose that $X+iY\in \ker (C+H)$, where $X,Y\in
\mathbb{R}^b$. We have
\begin{eqnarray}
(X^t-iY^t)\Big (R+iS\Big )(X+iY)= X^tRX+Y^tRY-2X^tSY\ge 0,
\label{eqn:ge0her}\end{eqnarray} for all $X,Y\in \mathbb{R}^b$. On
the other hand,
\begin{eqnarray}
\lefteqn{\Big ((A+R)+i(B+S)\Big )(X+iY)=}\nonumber \\
& & =\Big ((A+R)X-(B+S)Y\Big )+i\Big ( (B+S)X+(A+R)Y\Big ).
\label{eqn:0her}\end{eqnarray} Thus, $(A+R)X=(B+S)Y$ and
$(A+R)Y=-(B+S)X$. Let us multiply the first relation on the left
by $X^t$ and the second by $Y^t$, and sum: we get
$$\big (X^tRX+Y^tRY\big )+(X^tAX+Y^tAY)=2X^tSY.$$
In view of (\ref{eqn:ge0her}) and the positive definiteness of
$A$, we conclude $X^tAX=0$, $Y^tAY=0$ whence $X=Y=0$.

\bigskip

To complete the proof of Theorem \ref{thm:3}, we now only need to
observe that the asymptotic expansion described in the statement
holds for $\sum _{j=1}^{r_0}H_{3,j}^{(\omega,k)}$, in view of
Lemma \ref{lem:critical1} and the complex stationary phase Lemma.

\section{Proof of Theorem \ref{thm:2}.}

\noindent To ease the exposition, let us first prove Theorem
\ref{thm:2} under the following additional simplifying assumption:

\noindent \textit{$\omega \in \frak{h}^*_+$ lies in an elementary
fundamental wedge for the Hamiltonian action of $G$ on
$T^*X\setminus \{0\}$ (\cite{gs-hq}, page 357)}.

Recall that this means that $\omega$ is a regular value of the
moment map $\Psi :T^*X\setminus \{0\}\rightarrow \frak{g}^*$, and
that $G_\omega\subseteq G$ acts freely on $\Psi ^{-1}(\omega)$.

To begin with, let us then first of all clarify the relation
between the symplectic cone $Y$ and the submanifolds
$$W=:\Phi ^{-1}\left (C({\cal O})\right )\subseteq M,
\,\,Z=:\Psi ^{-1}\left (C({\cal O})\right )\subseteq T^*X\setminus
\{0\}.$$ In view of (\ref{eqn:momY}), if $p\in X$ and $r>0$ then
$$(p,r\alpha _p)\in Y\cap \Psi ^{-1}\left (C({\cal O})\right )\,
\Longleftrightarrow \,
\pi (p)\in \Phi ^{-1}\left ( C({\cal O})\right ).$$ In other
words, we have

\begin{lem} \label{lem:inters}
Let $\tilde \pi :Y\rightarrow M$ be the projection, $(x,r\alpha
_x)\mapsto \pi (x)$. Then
$$Y\cap Z\,=\,\tilde \pi ^{-1}\left (W \right ).$$\end{lem}

\noindent Furthermore, $0\not \in \Phi (M)$ clearly implies:
\begin{lem} \label{lem:proper}
$W\subseteq M$ is compact. \end{lem}

Recall, after \cite{gs-gq}, that $Z$ and $W$ are fibrating
coisotropic submanifolds of $T^*(X)\setminus \{0\}$ and $M$,
respectively. More precisely, there exist symplectic V-manifolds
in the sense of Satake
$$\left ( Z^\sharp,\Omega _{Z^\sharp} \right )
\mbox{ and }\left (W^\sharp ,\Omega _{W^\sharp} \right ),$$ and
smooth maps $p_Z:Z\rightarrow Z ^\sharp$, $p_W:W\rightarrow
W^\sharp$ whose fibres are the leaves of the corresponding null
foliations. Thus, if $\iota _Z:Z\rightarrow T^*(X)\setminus \{0\}$
and $\iota _W:W\rightarrow M$ are the inclusions, then $\iota
_Z^*(\Omega _{T^*X})=p_Z^*(\Omega _{Z^\sharp})$ and $\iota
_W^*(\Omega _M)=p^*_W(\Omega _{W^\sharp})$, where $\Omega
_{T^*X},\ldots$ denote the symplectic structures of the manifolds
$T^*X,\ldots$. The fibre of $p_Z$ (respectively, of $p_W$) through
a point $p\in Z$ (resp., $x\in W$) is the orbit of a certain
normal subgroup $H_f$ of the stabilizer subgroup $G_f$ under the
coadjoint action, where $f=\Psi (p)$ (resp., $f=\Phi (x)$). The
Lie algebra $\mathfrak{h}_f$ of $H_f$ is the codimension one ideal
of $\mathfrak{g}_f=\mathrm{Lie}(G_f)$ given by (\cite{gs-hq}, page
351): \begin{eqnarray}\mathfrak{h}_f=\left \{ \xi \in
\mathfrak{g}_f:<f,\xi>=0\right \}.\label{eqn:hf}\end{eqnarray}

Let $R$ be a symplectic manifold. Given an isotropic submanifold
$S\subset R\times R$, we shall denote by $S'$ the corresponding
isotropic relation, that is, the image of $S$ under the map
$(r_1,r_2)\mapsto (r_1,-r_2)$. Clearly $S'$ is an isotropic
submanifold of $R\times R^-$, where $R^-$ denotes $R$ with the
opposite symplectic structure. The fibre products
$$Z\times _{p_Z}Z\mbox{ and }W\times _{p_W}W$$
are then Lagrangian relations in $T^*(X\times X)$ and $T^*(M\times
M)$, respectively. We shall in a while be interested in the
composition of the conic isotropic relation $\Sigma '$, with
$\Sigma$ given by (\ref{eqn:wf}), with $D=:\left (Z\times
_{p_Z}Z\right )$. The previous discussion obviously implies:

\begin{lem}
$\Sigma '\circ D=\tilde \pi ^{-1}(W)\times _{p_Z}\tilde \pi
^{-1}(W)$. \label{lem:composition}\end{lem}

More explicitly, given that $\alpha$ is $G$-invariant we have
\begin{eqnarray}
\tilde \pi ^{-1}(W)\times _{p_Z}\tilde \pi ^{-1}(W)= \left \{\left
((x,r\alpha _x),(y,r\alpha _y)\right ): \pi
(x)\in W, \right .\nonumber \\
\left .y\in H_{\Phi (\pi (x))}\cdot x, r>0\right
\}.\label{eqn:circ}\end{eqnarray}

With the standard implicit identification between sections of
$L^{\otimes k}$ and $S^1$-equivariant funcions on $X$, let
$$H(X)_{(\omega)}=:\bigoplus _{k\ge 0}H^0(M,L^{\otimes
k})_{(\omega)}= \bigoplus _{k,\ell \ge 0}H^0(M,L^{\otimes
k})_{\ell\omega}\subseteq H(X)\subseteq L^2(X),$$ and let
$P_{(\omega)}:L^2(X)\rightarrow H(X)_{(\omega)}$ be the orthogonal
projector. If $\left \{s_j^{(k,\ell \omega)}\right \}$ is an
orthonormal basis of $H^0(M,L^{\otimes k})_{\ell \omega}$ for
$k,\ell =1,2,\ldots$, the Schwartz kernel of $P_{(\omega)}$ is
$$\tilde P_{(\omega)}(x,y)=\sum _{k,\ell}s_j^{(k,\ell \omega)}(x)
\otimes \overline
s_j^{(k,\ell \omega)}(y)\,\,\,\,\,\,(x,y\in X).$$ The k-th Fourier
component of $\tilde P_{(\omega)}$ is
$$\tilde P_{(\omega),k}(x,y)=\sum _\ell s_j^{(k,\ell
  \omega)}(x)\otimes
\overline s_j^{(k,\ell \omega)}(y) \,\,\,\,\,\,(x,y\in X),$$ and
we want to estimate asymptotically the diagonal behaviour of
$\tilde P_{(\omega),k}$ on $W=\Phi ^{-1}\left (C({\cal O})\right
)$.

To this end, we shall describe $P_{(\omega)}$ as a Fourier
integral operator with complex phase and study the geometry of the
associated (almost complex) canonical relation. In fact,
$P_{(\omega)}$ is the composition of two Fourier integral
operators whose associated canonical relations are related to the
$G$-action and to the complex structure, respectively. Namely, let
$$L^2(X)_{(\omega)}=\bigoplus _{\ell \ge 1}
L^2(X)_{\ell \omega}\subseteq L^2(X)$$
and denote by $Q_{(\omega)}:L^2(X)\rightarrow L^2(X)_{(\omega)}$
the orthogonal projector. If $\Pi$ denotes, as above, the Szeg\"o
projector, then $P_{(\omega)}=Q_{(\omega)}\circ \Pi$.

As mentioned already, $\Pi$ is an elliptic degree zero Fourier
integral operator with complex phase; its almost complex
Lagrangian submanifold $\tilde C$ is locally parametrized by the
phase function $\psi$, and can be characterized geometrically as
follows (\cite{bs}, Propositions 2.13 and 2.16). Let $\tilde X$
denote a complexification of $X$ \cite{ms}, so that the natural
(complex) symplectic structure on the cotangent bundle $T^*(\tilde
X)$ is the complexification of the symplectic structure on
$T^*(X)$. Let $\zeta =\{\zeta _j\}$ be the symbol of the boundary
Cauchy-Riemann operator $\overline \partial _b$, a smooth function
on $T^*(X)$, and let $\tilde \zeta=\{\tilde \zeta _j\}$ be its
almost analytic extension to $T^*(\tilde X)$.

\begin{thm} (\cite{bs})
Let $\Upsilon\subseteq T^*(\tilde X)\setminus \{0\}$ be the almost
analytic submanifold defined by $\tilde \zeta = 0$. Then
$\Upsilon$ is a conic involutive submanifold of $T^*(\tilde
X)\setminus \{0\}$, and $\tilde C$ is (up to almost analytic
equivalence) the unique Lagrangian submanifold of $\left
(T^*(\tilde X)\setminus \{0\}\right )\times \left (T^*(\tilde
X)\setminus \{0\}\right )$ contained in $\Upsilon\times \overline
\Upsilon$ and containing the wave front (\ref{eqn:wf}) of
$\Pi$.\label{thm:upsilon}\end{thm}

Next, $Q_{(\omega)}$ is also an elliptic degree zero Fourier
integral operator, with real phase however. Its associated
canonical relation is described in \cite{gs-hq}, Theorem 6.7, in
terms of the moment map $\Psi$. Namely, with $\omega$ and $C({\cal
O})$ as in Definition \ref{defn:cone}, assume that $\Psi$ is
transversal to $C({\cal O})$. The conic Lagrangian relation
$$D\subseteq \left (
  T^*(X)\setminus \{0\}\right )\times
 \left (T^*(X)\setminus \{0\}\right )$$ associated to $Q_{(\omega)}$ is then the
fibre product $Z\times _p Z$. To discuss the composition of $\Pi$
and $Q_{(\omega)}$, we shall view the latter as a Fourier integral
operator with complex phase, whose associated almost complex
canonical relation $$\tilde D\subseteq \left (T^*(\tilde X)\times
\{0\}\right )\times \left (T^*(\tilde
  X)\setminus \{0\}\right )$$ is simply the almost analytic extension of $D$.
In the language and notation of almost analytic machinery
\cite{ms}, this can be described as follows.

The moment map $\Psi :T^*(X)\rightarrow \frak{g}^*$ extends almost
analytically to $\tilde \Psi :T^*(\tilde X)\rightarrow
\frak{g}_c^*$, where $\frak{g}_c=\frak{g}\otimes
_{\mathbb{R}}\mathbb{C}$. Let $\tilde C(\tilde {\cal O})\subseteq
\frak{g}^*_c$ be the complexification of $C({\cal O})\subseteq
\frak{g}^*$; this can be explicitly decribed as
\begin{eqnarray*}
\tilde C(\tilde {\cal O})=\left \{\lambda \omega':\lambda \in
  \mathbb{C}^*,\,
\omega'\in \tilde {\cal O} \right \},\end{eqnarray*} where $\tilde
{\cal O}\subseteq \frak{g}^*_c$ is the coadjoint orbit of $\omega$
in $\frak{g}_c^*$ under the complexification $\tilde G$ of $G$.
Then $\tilde \Psi$ is transversal to $\tilde {\cal O}$, and the
almost analytic submanifold $\tilde Z=\tilde \Psi ^{-1}(\tilde
{\cal O})\subseteq T^*(\tilde X)\setminus \{0\}$ is the
complexification of $Z$. The fibration $p:Z\rightarrow  M_\omega
^\sharp$ extends almost analytically to $\tilde p:\tilde
Z\rightarrow \tilde M_\omega ^\sharp$, and the fibre product
$$\tilde D=:\tilde Z\times _{\tilde p}\tilde Z\subseteq T^*(\tilde
X)\times T^*(\tilde X)$$ is the almost analytic extension of
$D=Z\times _pZ$. It is therefore a Lagrangian relation in
$T^*(\tilde X)\times T^*(\tilde X)$.

Clearly, $D$ and $\Sigma$ are the sets of real points of $\tilde
D$ and $\tilde C$, respectively. Set $$\Delta =T^*( X)\times {\rm
diag}\left(T^*( X)\right )\times T^*( X)$$ and
$$\tilde \Delta=T^*(\tilde X)\times {\rm diag}\left (T^*(\tilde
  X)\right )\times
T^*(\tilde X).$$

\begin{lem} \label{lem:transv}
$\tilde D\times \tilde C$ and $\tilde \Delta$ meet transversally
in $T^*(\tilde X)\times T^*(\tilde X)\times T^*(\tilde X)\times
T^*(\tilde X)$ at every point of $\left (D \times \Sigma\right
)\cap \Delta$.\end{lem}

\noindent {\it Proof.} All statements in the following arguments
are meant to be local along the real locus. We have
\begin{eqnarray*}
\tilde D=\left \{(w,w')\in \left (T^*(\tilde X)\setminus
\{0\}\right )
\times \left (T^*(\tilde X)\setminus \{0\}\right ): \right. \nonumber \\
\Psi (w),\Psi (w')\in \tilde C( \tilde {\cal O}), \left .\tilde p
(w)=\tilde p(w')\right \}.\end{eqnarray*} If $(w,w')\in \tilde D$,
then the tangent space to $\tilde D$ at $(w,w')$ is
\begin{eqnarray*}
T_{(w,w')}(\tilde D)=\left \{(v,v')\in T_{(w,w')}(\tilde Z\times
\tilde Z ):\, d_w\tilde p(v)=d_{w'}\tilde p(v')\right
\}.\end{eqnarray*} With $g=\dim (G)$, set $k=g-\dim ({\cal
O})-1={\rm codim}\big (C({\cal O})\big )$. If $v'\in T_{w'}(\tilde
Z)$, the collection of all $v\in T_w(\tilde Z)$ with $d_w\tilde
p(v)=d_{w'}\tilde p(v')$ is an affine space of dimension $k$.

Suppose next $y=(w,w',w',w'')\in \left (D \times \Sigma\right
)\cap \Delta$. Then
\begin{eqnarray*}
\left (T_{(w,w')}\big (\tilde D\big )\times T_{(w',w'')}\big
(\tilde
  C\big )\right )\cap T_y\big (\tilde \Delta \big )=\left
  \{(v,v',v',v''):(v,v')
\in T_{(w,w')}(\tilde Z\times
\tilde Z ), \right .\\
\left .d_w\tilde p(v)=d_{w'}\tilde p(v'),\,(v',v'')\in
T_{(w',w'')}(\tilde C)\right \}.
\end{eqnarray*}
The condition $(v,v')\in T_{(w,w')}(\tilde Z\times \tilde Z )$ may
be rewritten
$$d_w\tilde \Psi (v)\in T_{\tilde \Psi (w)}\big (\tilde C(\tilde {\cal
  O})\big ),
\,\,\,d_{w'}\tilde \Psi (v')\in T_{\tilde \Psi (w')}\big (\tilde
C(\tilde {\cal O})\big ).$$
\begin{claim}
$\tilde C$ is transversal to $\tilde C(\tilde {\cal O})$ under the
map $\tilde \gamma: (w,w')\mapsto \tilde \Psi (w)$ (along the real
locus).\end{claim}

\noindent {\it Proof.} Certainly $\tilde C\supseteq \tilde
\Sigma\supseteq \Sigma$, where
$$
\tilde \Sigma =\left \{(p,\lambda \tilde \alpha _p,p,-\lambda
\tilde \alpha _p): p\in \tilde X,\,\lambda \in \mathbb{C}^* \right
\}$$ is the complexification of $\Sigma$; here $\tilde \alpha$ is
the complexification of the connection 1-form $\alpha$. It
suffices to show that $\tilde \Sigma$ is transversal to $\tilde
C(\tilde {\cal O})$ under $\tilde \gamma$, whence that $\Sigma$ is
transversal to $C({\cal O})$ under the map $\gamma :(w,w')\mapsto
\Psi (w)$. By homogeneity of the moment map on the cotangent
bundle and equality (A.11) in the Appendix to \cite{gs-gq}, if
$p\in X$ and $r>0$ then
\begin{eqnarray*}
\gamma \left ((p,r \alpha _p,p,-r \alpha _p) \right )= r\Psi \left
((p, \alpha _p)\right ) =\Phi (p).\end{eqnarray*} The statement
follows from the hypothesis that $\Phi$ be transversal to $C({\cal
O})$.

\bigskip

Returning to the proof of Lemma \ref{lem:transv}, the Claim
clearly implies that
$$\dim _{\mathbb{C}}\left (\left (T_{(w,w')}\big (\tilde D\big )\times
    T_{(w',w'')}
\big (\tilde C\big )\right )\cap T_y\big (\tilde \Delta \big
)\right )=\dim (\tilde C)=2r,$$ where $r=2n+1=\dim(X)$. We have on
the other hand
\begin{eqnarray*} \dim _{\mathbb{C}}(\tilde C)+
\dim _{\mathbb{C}}(\tilde D)- {\rm codim}_{\mathbb{C}}(\tilde
\Delta)=2r+2r-2r=2r,
\end{eqnarray*} and thus the statement.

\begin{lem} \label{lem:clean} The projection $\Pi=(\pi _1,\pi _4)$
onto the the first and fourth factors,
$$\Pi:(w,w',w',w'')\in (D\times \Sigma)\cap \Delta
\mapsto (w,w'')\in \left (T^*(X)\setminus
  \{0\}\right )\times \left (T^*(X)\setminus \{0\}\right ),$$
is injective and proper.
\end{lem}
 \noindent {\it Proof.} Injectivity follows immediately from the
 description of $\Sigma$.

Let next $K\subset \left (T^*(X)\setminus \{0\}\right ) \times
\left  (T^*(X)\setminus \{0\}\right )$ be a compact subset. We
want to establish that $\Pi ^{-1}(K)$ is compact. Now
\begin{eqnarray*}
(D\times \Sigma)\cap \Delta =\left \{\left (h\cdot (p,r\alpha _p),
(p,r\alpha _p),(p,r\alpha _p),(p,-r\alpha _p)\right ):\right .
\\ \left .\pi (p)\in
W, r>0,h\in H_{\Phi \circ \pi (p)}\right \}.\end{eqnarray*} Here
$H_f$ denotes the codimension one closed Lie subgroup of the
stabilizer $G_f\subseteq G$ of an element $f\in \mathfrak{g}^*$
discussed on page 349 of \cite{gs-hq}. The projection
$K'=q_2(K)\subseteq T^*(X)\setminus \{0\}$ of $K$ onto the second
factor is a compact subset of $T^*(X)\setminus \{0\}$. Clearly,
$\Pi ^{-1}(K)$ is a closed subset of
\begin{eqnarray}
\label{eqn:p-1K} \left \{\left (h\cdot (p,r\alpha _p),(p,r\alpha
_p),(p,r\alpha _p),(p,-r\alpha _p)\right ):\right .\pi (p)\in W,
r>0,\\ \left . h\in H_{\Phi \circ \pi (p)},(p,-r\alpha _p)\in
K'\right \}.
\end{eqnarray}
It thus suffices to show that the latter set is compact. It is
obvious that the union $\bigcup _{\pi (p)\in W}H_{\Phi \circ \pi
(p)}$ is a compact subset of $G$. Since $W$ is compact, so is $\pi
^{-1}(W)\subseteq X$. Therefore the set of all $r>0$ such that
$(p,-r\alpha _p)\in K'$ for some $p\in \pi ^{-1}(W)$ is contained
in  a closed interval $[a,b]\subset \mathbb{R}_{+}$. The statement
follows.

\begin{rem} \label{rem:compact} With a view to Corollary \ref{cor:M'},
we remark that
the condition that $W$ is compact is not essential in the proof of
Lemma \ref{lem:clean}. In fact, the projection $K^{\prime
\prime}\subseteq M$ of $K'$ under the composition $T^*X\rightarrow
X\stackrel{\pi}{\rightarrow} M$ is at any rate a compact subset of
$W$, and the condition $\pi (p)\in W$ may be replaced by the
condition $\pi (p)\in K^{\prime \prime}$ in
(\ref{eqn:p-1K}).\end{rem}

\noindent Given this, we can compose the Fourier integral
operators $Q_{(\omega)}$ and $\Pi$ \cite{ms}:

\begin{cor} $P_{(\omega)}$ is a degree zero elliptic Fourier
integral operator with complex phase, associated to the almost
analytic Lagrangian submanifold
$$\Gamma =\left (\tilde D\circ \tilde C'\right )'\subseteq
\left(T^*(\tilde X) \setminus \{0\} \right )\times \left
(T^*(\tilde X)\setminus \{0\}\right ).$$\end{cor}

We now characterize $\Gamma$ geometrically. Let $\Upsilon
\subseteq T^*(\tilde X)\setminus \{0\}$ be as in the statement of
Theorem \ref{thm:upsilon}.

\begin{lem} Up to almost analytic equivalence, $\Gamma$ is uniquely determined by
the properties: \label{lem:unique} \noindent i) $\Gamma \subseteq
\left (T^*(\tilde X)\setminus \{0\}\right )\times \left
(T^*(\tilde X)\setminus \{0\}\right )$ is an almost analytic
Lagrangian submanifold;

\noindent ii) $\Gamma '\supseteq D\circ \Sigma '$;

\noindent iii) $\Gamma \subseteq \Upsilon \times \overline
\Upsilon$.\end{lem}

\noindent {\it Proof.} It is clear that $\Gamma '\supseteq D\circ
\Sigma$. Let us show that $\Gamma \subseteq \Upsilon \times
\overline \Upsilon$. We may equivalently show that $\Gamma
'=\tilde D\circ \tilde C' \subseteq \Upsilon \times \overline
\Upsilon$. We have, by definition,
\begin{eqnarray*}
\Gamma '=\left \{ (w,w')\in \left (T^*(\tilde X)\setminus
\{0\}\right )\times \left (T^*(\tilde X)\setminus \{0\}\right ):
\exists \, w''\in T^*(\tilde X)\right .\\
\left .\mbox{ such that }(w,w'')\in \tilde D,\, (w'',w')\in \tilde
C\right \}.
\end{eqnarray*}
Since $\tilde C\subseteq \Upsilon \times \overline \Upsilon$, we
have $w'\in \overline \Upsilon$ for all $(w,w')\in \Gamma$. To
show that $w\in \Upsilon$, recall that
$$\tilde D=\tilde Z\times _{\tilde p}\tilde Z\subseteq T^*(\tilde
X)\times T^*(\tilde X).$$ It is thus sufficient to show that if
$w,w''\in \tilde Z$, $\tilde d(w'')=0$ and $\tilde p (w)=\tilde
p(w'')$, then $\tilde d(w)=0$. Being the complexification of the
isotropic submanifold $Z$, $\tilde Z$ is also isotropic, and the
leaf of its null foliation through $p\in \tilde Z$ is the orbit of
the complex group $\tilde H_f\subseteq \tilde G$. Since the action
of $G$ preserves the horizontal distribution $H(X/M)$
(equivalently, the connection 1-form $\alpha$), the symbol $d$ of
the boudary Cauchy-Riemann operator is $G$-invariant.
Complexifying, the action of $\tilde G$ preserves the horizontal
distribution $H(\tilde X/\tilde M)$, whence the components $\tilde
d_j$ of $\tilde d$ are $\tilde G$-invariant. Equivalently, since
$d$ is $G$-invariant the Hamiltonian vector fields $X_{d_j}$ are
in the kernel of $dp$; their complexifications $\tilde X_{\tilde
d_j}$'s are therefore in the kernel of $d\tilde p$.

To establish the Lemma, we are thus reduced to proving that if
$$\Lambda \subset \left (T^*(\tilde X)\setminus \{0\}\right )\times
\left (T^*(\tilde X)\setminus \{0\}\right )$$ is an almost
analytic Lagrangian submanifold such that $$\Lambda '\supseteq
D\circ \Sigma ' \mbox{ and }\Lambda \subseteq \Upsilon \times
\overline \Upsilon,$$ then (up to almost analytic equivalence)
$\Gamma =\Lambda$.

Clearly $\Lambda '\supseteq (D\circ \Sigma ')\tilde{}$, where
$(D\circ \Sigma ')\tilde{}$ denotes the complexification of
$D\circ \Sigma '$. Now $\left ((D\circ \Sigma ')'\right )\tilde{}$
is a $2(n+1)$-dimensional almost analytic isotropic submanifold of
$\left (T^*(\tilde X)\setminus \{0\}\right )\times \left
(T^*(\tilde X)\setminus \{0\}\right )$. On the other hand,
$\Upsilon \times \overline \Upsilon$ is the coisotropic
submanifold defined by the $2n$ equations $\tilde d_j^{(1)}=0$ and
$\overline {\tilde d_j}^{(2)}=0$, $j=1,\ldots,n$ ($(1)$ and $(2)$
stand for the component of evaluation). The leaves of the null
foliation of $\Upsilon\times \overline \Upsilon$ are then
generated by the flow of the Hamiltonian vector fields $\tilde
X_{\tilde d_j}^{(1)}$'s and $\tilde X_{\overline {\tilde
d}_j}^{(2)}$'s. These span at each point of $\left ((D\circ \Sigma
')'\right )\tilde{}$ a $2n$-dimensional linear space which is
trasversal to the tangent space of $(D\circ \Sigma ')\tilde{}$.
Hence $\Lambda$ must be the union of the leaves passing through
$\left
  ((D\circ \Sigma ')'
\right )\tilde{}$.

\bigskip
We now produce suitable local holomorphic coordinates on $M$ in
the neighbourhood of $p\in W$, and on $L$ in the neighbourhood of
$x\in Z$. These coordinates are adapted to the fibrating structure
of $W$ and $Z$.

\begin{lem} For any $p\in W$ and any sufficiently small open
neighbourhood $U\ni x$ in $W$, the image of $U$ under $p_W$ is an
open K\"ahler manifold $S\subseteq W^\sharp$. The map $p_U=\left
.p\right |_U:U\rightarrow S$ is a submersion, and $\left .\Omega
\right |_U=p_W^*\left (\Omega _S\right )$, where $\Omega _S$ is
the K\"ahler structure of $S$. \label{lem:image} \end{lem}

\noindent {\it Proof.} This is proved by a straightforward
adaptation of arguments from \cite{gs-gq} and \cite{gs-hq}; we
sketch the proof for the reader's convenience and to introduce
some notation. As already recalled, by the theory in \cite{gs-hq},
we have the following picture: If $p\in W$ and $f=\Phi (p)\in
\frak{g}^*$, let $H_f\subseteq G$ be the normal closed Lie
subgroup of the stabilizer $G_f$ of $f$ with Lie algebra
(\ref{eqn:hf}). Then the fibre $p_W^{-1}(p)\subseteq W$ is the
orbit of $p$ under $H_f$. Therefore, $p^{-1}_W(p)$ is
diffeomorphic to the quotient $H_f/H(p)$, where $H(p)\subseteq
H_f$ is the stabilizer of $p$ in $H_f$. Furthermore,
$H(p)\subseteq H_f$ is a discrete subgroup of $H_f$. Thus, $T_pW$
is an coisotropic subspace of $(T_pM,\Omega _p)$, with symplectic
complement $\left (T_pW\right )^\perp =T_p(H_f\cdot p)\subseteq
T_pW$. Let us set $R_p=:T_p(H_f\cdot p)$.

We conclude that $p_W$ is a submersion and $S$ is a symplectic
manifold of real dimension $2\left (n-\dim (H_f)\right )=2\left
(n-g+\dim (\mathcal{O})+1\right )$. Furthermore, for $p\in M$ let
$F_p\subset T_pM\otimes \mathbb{C}$ be the $+i$-eigenspace of the
complex structure $J_p\in {\rm End}\left (T_pM\right )$. By Lemma
3.6 of \cite{gs-gq},
\begin{equation}F_p\cap \left (R_p\otimes \mathbb{C}\right )=0.
\label{eqn:transv}\end{equation}

For $p\in W$, we set
$$F'_p=:F_p\cap \left (T_pW\otimes \mathbb{C}\right );$$
then $(F'_p)^\perp$ has complex dimension $n+\dim H_f$ ($n$ is the
complex dimension of $M$). Therefore, by $H_f$-invariance and
dimension count, $F'$ descends to a well-defined positive
Lagrangian complex distribution on $S$. Integrability follows as
in {\it loc. cit.}, and therefore $S$ inherits compatible complex
and symplectic structures.

\begin{lem} Pick any $p\in W$ and let $U$ be any sufficiently small open
neighbourhood of $p$ in $W$. Set $S=p_W(U)$, and endow it with the
holomorphic structure described in Lemma \ref{lem:image}. Then
there exists a map $\sigma :S\rightarrow W$ satisfying the
following properties: i) $\sigma \left (p_W(p)\right )=p$; ii):
$\sigma$ is holomorphic as a map from $S$ to $M$; iii): $\sigma$
is a section of $p_W$, that is, $p_W\circ \sigma ={\rm id}_{S}$.
\label{lem:section}\end{lem}

\noindent {\it Proof.} The only property not immediately obvious
is (ii). Set as above $k=g-\dim ({\cal O})-1$; we may choose local
coordinates $t_1,\ldots, t_g$ for $\mathfrak{g}^*$ on a
neighbourhood $V$ of $f=\Phi (p)$ centered at $f$ and such that
$C({\cal O})\cap V$ is defined by $t_1=\cdots =t_k=0$. Let
$U\subseteq \Phi ^{-1}(V)$ be a sufficiently small neighbourhood
of $p$. The map $$\psi :U\rightarrow \mathbb{R}^k,\, p'\in
U\mapsto \left (t_1(\Phi (p)),\ldots,t_k(\Phi (p)\right )\in
\mathbb{R}^k$$ is a submersion, and its image a neighbourhood $D$
of $0\in \mathbb{R}^k$. For $c\in U$, let $W_c=\psi
^{-1}(c)\subset U$. Then $W_0=W\cap U$, and $W_c$ is a submanifold
of $U$ for every $c\in D$. If $q\in U$, define
$$F'_q=:F_q\cap \left (T_q(W_{\psi (q)})\otimes \mathbb{C}\right ).$$
Let us set $R_q=:T_q(W_{\psi (q)})^\perp$. Given
(\ref{eqn:transv}), we have
\begin{equation}F_q\cap \left (R_q\otimes \mathbb{C}\right )=0
\label{eqn:transv1}\end{equation} for every $q\in U$ (perhaps
after further restricting $U$). This implies that $F'$ is a
complex distribution on $U$, of rank $n-g+\dim (\mathcal{O})+1$.
It is integrable, essentially by definition, and therefore the
complex version of the Frobenius integrability theorem applies.
Thus, we may find local holomorphic coordinates $r_1,\ldots,r_n$
on $U$ centered at $p$, such that $F'={\rm span}\left
\{\partial/\partial r_j\,: \, j=1,\ldots,k\right \}$. After
further restricting if necessary, the manifold $\tilde
S=\{r_{k+1}=\cdots=r_n=0\}$ is then a section of $p$, as required.

\bigskip

The $G$-invariance of $L$ and $\nabla$ clearly implies

\begin{lem} Let $U\subseteq Z$ and $S=p(U)\subseteq M_\omega ^\sharp$
be as in Lemma \ref{lem:image}. On $S$ there are an hermitian line
bundle $\big (L_S,h_{L_S}\big )$ with (unique) compatible
connection $\nabla _{L_S}$, such that $\left .L\right
|_U=p_U^*(L_S)$, $h_{L_S}=p_U^*(h_L)$ and $\nabla _L=p_U^*\left
(\nabla _{L_S}\right )$. \label{lem:pullback}\end{lem}

\noindent We have in fact $L_S=\sigma ^*(L)$, $h_{L_S}=\sigma
^*(h_L)$, $\nabla _{L_S}=\sigma ^*\left (\nabla _{L_S}\right )$.

Let us now fix $p\in W$ and let $U$ be the open neighbourhood of
$p$ in $W$, $S$ the open K\"ahler manifold and $\pi _U=\left
.p_W\right |_U:U\rightarrow S$ be the submersion described in
Lemma \ref{lem:image}. Let $\sigma :S\rightarrow U\subseteq M$ be
the holomorphic section described in Lemma \ref{lem:section}. Set
$\overline p=\pi _U(p)\in S$, and let $(z_1,\ldots,z_k)$ be local
holomorphic coordinates on $S$ centered at $\overline p$, where
$k=n-g+\dim {\cal O}+1$. Let $H_f^c\subseteq \tilde G$ be the
complexification of $H_f$, with Lie algebra
$\frak{h}_f^c=\frak{h}_f\otimes \mathbb{C}$. Let $\exp
:\frak{h}_f^c\rightarrow H_f^c$ be the exponential map of $H_f^c$,
and let $V\subset \frak{h}_f^c$ be an open neighbourhood of the
origin such that $\exp$ induces a diffeomorphism of $V$ onto its
image. Then, perhaps retricting $S$ and $V$ to smaller open
subsets, the holomorphic map $\zeta :S\times V\rightarrow M$ given
by
$$\zeta: (s,\upsilon)\mapsto \exp (\upsilon)\cdot \sigma (s)$$
is a diffeomorphism onto its image. We obtain local holomorphic
coordinates $(z_1,\ldots,z_k,w_1,\ldots,w_h)$ on $M$ centered at
$x$; here $w_j=u_j+iv_j$, where $u_j$ and $v_j$ are linear
coordinates on $\frak{h}_f$, and $h=n-k=g-\dim {\cal O}-1$.

We may also suppose without loss of generality that the local
holomorphic coordinates $(z,w)$ are preferred at $p$ in the sense
of \cite{sz}: in our integrable case this simply means that in the
given coordinate system the K\"ahler form at $p$ is the standard
symplectic structure $\omega _0$ on $\mathbb{C}^n$. To see this,
notice that we have a symplectically orthogonal direct sum
\begin{eqnarray*}
T_p(M)\otimes \mathbb{C}\,=\, F'_p\oplus \left [R_p\oplus
J_p(R_p)\right ],\end{eqnarray*} where $R_p=T_p(H_f\cdot p)$.
Furthermore, $R_p$ and $J_p(R_p)$ are dually paired Lagrangian
subspaces of the symplectic subspace $R_p\oplus J_p(R_p)$. In
terms of the decomposition
$$T_0(S\times V)=\mathbb{C}^k\oplus \frak{h}_f^c\,\mbox{ and
}\,\frak{h}_f^c= \frak{h}_f\oplus i\cdot \frak{h}_f$$ we have
\begin{eqnarray*}F_p'=d_0\zeta \left (\mathbb{C}^k\oplus \{0\}\right ),\,R_p=
d_0\zeta \left (\{0\}\oplus \frak{h}_f\oplus \{0\}\right ),\\
J_p(R_p)=d_0\zeta \left (\{0\}\oplus \{0\}\oplus i\cdot
\frak{h}_f\right ).
\end{eqnarray*}
Thus, in order for $(z,w)$ to be a preferred system of local
holomorphic coordinates at $p$ it is sufficient to choose the
$z_i$'s so that they form a system of preferred coordinates on $S$
at $\overline p$, and to make an appropriate choice for a basis in
$\frak{h}_f$.

Let now $e_L$ be a preferred holomorphic local frame for $L$ at
$p$ in the sense of \cite{sz}. Given our choice of a preferred
system of local holomorphic coordinates $(z,w)$ at $p$, the
function $\beta =||e_L^*||^2=||e_L||^{-2}$ then satisifes
\begin{eqnarray}\beta (z,w)=1+|z|^2+|w|^2+\cdots.
\label{eqn:expansion}\end{eqnarray} A point in the neighbourhood
of $L^*(p)\subseteq L^*$ with local holomorphic coordinates
$(z,w,\lambda)$ lies on $X$ if and only if $|\lambda |^2 \beta
(z,w)=1$.

In order to obtain an asymtpotic expansion for the equivariant
Sz\"ego kernels $P_\omega$ on the diagonal, we shall need a more
complete description of their microlocal structure near any given
point of ${\rm diag}(W)$. To this end, we shall now prove the
existence of a regular positive phase function, in the sense of
\cite{ms}, microlocally associated to the Lagrangian almost
analytic manifold $\Gamma$ in the neighbourhood of any point in
$(p,p)\in \mathrm{diag}(W)$. The phase function will be (locally)
defined on $X\times X\times \mathbb{R}^+$ and have the form
$t\cdot \psi _\omega (x,y)$.

Let us then fix $p_0\in W$ and adopt in an open neighbourhood $T$
of $p_0$ in $M$ the local holomorphic coordinates $(z,w)$
introduced above. This determines local holomorphic coordinates
$(z,w,z',w')$ on the open neighbourhood $T\times T\subseteq
M\times M$ of $(p_0,p_0)\in \mathrm{diag}(W)$. Let us choose on
$T$ (perhaps after restriction) a preferred holomorphic local
frame $e_L$ for $L$ at $p_0$.

Recall that $\tilde M$ and $\tilde X$ denote the almost analytic
extensions of $M$ and $X$, and that the $S^1$-principal bundle
$\pi :X\rightarrow M$ almost analytically extends to a
$\mathbb{C}^*$-bundle $\tilde \pi :\tilde X\rightarrow \tilde M$.

\begin{prop}\label{prop:phasegamma}
Let $p_0\in W$ and $T$ be as above. Then, perhaps after
restricting $T$ to a smaller open neighbourhood of $p_0$, there
exists a smooth function $b\in {\cal C}^\infty (T\times T)$ with
the following property: Let $\psi \in {\cal C}^\infty (\left
.X\right |_T\times \left .X\right |_T)$ be the restriction of
$\hat \psi \in {\cal C}^\infty (\left .L\right |_T\times \left
.L\right |_T)$ given by:
$$\hat \psi (\ell,\ell ')=i\cdot
\left (1-\overline \lambda \mu \,  b(p,q)\right ),$$ where
$\ell,\,\ell '\in \left .L\right |_T$ correspond to
$(p,\lambda),\, (q,\mu)\in T\times \mathbb{C}$, respectively. Then
$t\cdot \psi \in {\cal C}^\infty \left (\left . X\right |_T\times
\left .X\right |_T \times \mathbb{R}_+\right )$ is a regular phase
function for $\Gamma$.
\end{prop}

This means the following \cite{ms}, \cite{sz}. Let $\tilde \psi$
be the almost analytic extension of $\psi$ to $\left .\tilde X
\right |_{\tilde T} \times \left .\tilde X\right |_{\tilde T}$.
Then
\begin{eqnarray}
\Gamma \,=\, \left \{(\tilde x,\, td_{\tilde x}\tilde \psi,\tilde
y,td_{\tilde y}\tilde \psi)\,: \,\tilde \psi (\tilde x,\tilde
y)=0,\, t\in \mathbb{R}_+\right \}.
\label{eqn:phasefunction}\end{eqnarray} In addition, we shall
prove in Proposition \ref{prop:positivetype} that $\psi$ is of
positive type, that is, the imaginary part of $\psi$ satisfies $
{\cal I}(\psi)\ge 0$.

\noindent \textit{Proof of Proposition \ref{prop:phasegamma}.} We
shall adapt the reduction argument in the proof of Theorem 2.1 of
\cite{sz} to the present equivariant context.

As a first step, we shall replace the holomorphic preferred local
section $e_L$ of $L$ at $p$ with a different local section,
$\tilde e_L$, not holomorphic but capturing the equivariant
geometry of our picture. We shall construct a phase function for
$\Gamma$ of the asserted form $\phi =i\left (1 -\lambda \overline
\mu c(p,q)\right )$ in the local (nonholomorphic) coordinates on
$L^*$ (and $X$) associated to the holomorphic coordinates $(z,w)$
on $M$ and the section $\tilde e_L^*$. The phase function in the
original holomorphic coordinates will then be given by pull-back
under the change of basis from $e_L^*$ to $\tilde e_L^*$. Going
back to holomorphic coordinates will be convenient in the proof of
Proposition \ref{prop:positivetype}.

To this end, in the notation of Lemmas \ref{lem:section} and
\ref{lem:image}, let us first of all choose a local holomorphic
section $e_{L_S}$ of the Hermitian holomorphic line bundle
$L_S=\sigma ^*\left (L\right )$ on $S$ which is preferred at
$p_0$. We shall view this as a section of $L$ defined on the
submanifold $\sigma (S)\subset M$.

Next we shall extend $e_{L_S}$ to a local section of $L$ by using
the action of the complexified groups $\tilde H_f$. More
precisely, set
\begin{eqnarray*}
\mathcal{\tilde H}=\left \{(h,s)\in \tilde G\times S\,  :\, h\in
\tilde H_{\Phi (\sigma (s))}\right \},\end{eqnarray*}
\begin{eqnarray*}
\mathcal{H}=\left \{(h,s)\in G\times S\, :\, h\in H_{\Phi (\sigma
(s))}\right \}. \end{eqnarray*} Then $\mathcal{H} \subseteqq
\mathcal{\tilde H}$ and the projection onto the second factor
induces smooth fibrations $\tilde r:\mathcal{\tilde H}\rightarrow
S$ and $r:\mathcal{H}\rightarrow S$ with $r^{-1}(s)=H_{\Phi
(\sigma (s))}$ and $\tilde r^{-1}(s)=\tilde H_{\Phi (\sigma (s))}$
for every $s\in S$. There is an obvious smooth map
$$A:\tilde {\cal H}\longrightarrow M,\,\,\,\,\,
(h,s)\mapsto h\cdot \sigma (s)$$ which is a local diffeomorphism
along the unit section ${\cal E}=\{(e,s):s\in S\}$ of $\tilde r$.
Therefore, $A$ induces a diffeomorphism of an open neighbourhood
$\tilde {\cal U}$ of $\cal E$ onto an open neighbourhood $T$ of
$\sigma (S)$ in $M$. Let ${\cal U}=\tilde {\cal U}\cap {\cal H}$.
Then $\cal U$ maps diffeomorphically under $A$ onto an open
neighbourhood $T'=T\cap W$ of $\sigma (S)$ in $W$.

$A$ is covered by the smooth map
$$B:\tilde {\cal H}\longrightarrow L,\,\,\,\,\,(h,s)
\mapsto h\cdot e_{L_S}\left (\sigma (s)\right ).$$ Setting $\tilde
e_{L}\left (A(h,s)\right )=B(h,s)$ ($(h,s)\in \mathcal{\tilde H}$)
we thus define a $\mathcal{C}^{\infty}$ section of $L$ over $T$.

Let $\tilde e_{L_{S}}^{*}$ be the dual frame of $L^*$ over $T$,
and set $\tilde \beta =\left \| \tilde e_{L_{S}}^{*}\right \|
^2=\left \| \tilde e_{L_{S}}\right \| ^{-2}$. Since $G$ preserves
the Hermitian structure of $L$, $\tilde \beta$ is constant along
the fibres of $p_W:W\rightarrow W^{\sharp}$ in $T$. Then $\tilde
\beta ^{-1/2}\tilde e_{L_{S}}^{*}$ is a unitary section of $L^*$
over $T$, and induces a trivialization $\chi : S^1\times T\cong
\left .X\right |_T$,
$$\chi (\lambda ,p)\, = \,
\left (p,\frac{\lambda}{\sqrt{\tilde \beta (p)}}\tilde
e_{L_{S}}^{*}(p)\right ).$$

Suppose that $p'\in T$ has local coordinates $(z,a+i0)$ ($a\in
\mathbb{R}$). Pick $x'\in \pi ^{-1}(p')$ with coordinates
$(z,a,\lambda)$. If $h\in H_{\Phi (p')}$ and $h\cdot p'\in T$,
then $h\cdot p'\in T$ and it has coordinates $(z,b,\lambda)$ for
some $b\in \mathbb{R}$.

Now, in view of Lemma \ref{lem:unique} and (\ref{eqn:circ}) $\psi$
will parametrize $\Gamma$ if the following conditions are
satisfied:

i) set $Z(\psi)=\left \{(x,y)\in \left .X\right |_{T}\times \left
. X\right |_{T}: \psi (x,y)=0\right \}$; then $$Z(\psi)=\left
\{(x,h\cdot x)\in \left .X\right |_{T}\times \left .X\right |_{T}:
\pi (x)\in W, h\in H_{\Phi (\pi (x))}\right \};$$

ii): $\left .d_x\psi \right |_{Z(\psi)}=q_1^*(\alpha)$, $-\left
.d_y\psi \right |_{Z(\psi)}=q_2^*(\alpha)$, where $q_i:X\times
X\rightarrow X$ is the projection onto the $i$-th factor;

iii): passing to almost analytic extensions (and writing $\tilde
X$ for $\left .\tilde X\right |_{\tilde T}$), set $\tilde Z(\tilde
\psi)= \left \{(\tilde x,\tilde y)\in \tilde X\times \tilde X:
\tilde \psi (\tilde x,\tilde y)=0\right \}$; then $\tilde \zeta
\left (\tilde x,d_{\tilde x}\tilde \psi \right )=0$, $\tilde \zeta
\left (\tilde y,d_{\tilde y}\tilde \psi \right )=0$ on $\tilde
Z(\tilde \psi)$; here $\zeta$ is the symbol of boundary $\overline
\partial$-operator on $X$, $\tilde \zeta$ is the almost analytic
extension of $\zeta$ (and, as in Theorem \ref{thm:upsilon} and
Lemma \ref{lem:unique}, $\Upsilon$ is defined by $\tilde \zeta
=0$).

We look for a solution of the form $\phi =i\left (1-\lambda
\overline \mu c(p,q)\right )$, and view i), ii) and iii) as
conditions on $c$. To determine $c$, we shall symplectically
reduce $\Gamma$ with respect to the $S^1$-simmetry.

Namely, the action of $S^1$ on $X$ lifts to a Hamiltonian action
on $T^*(X)\setminus \{0\}$. Let $\nu :T^*X\setminus
\{0\}\rightarrow \mathbb{R}$ be the moment map. These actions in
turn almost analytically extend to actions of $\mathbb{C}^*$ on
$\tilde X$ and $T^*(\tilde X)$, and $\nu$ to an almost analytic
map $\tilde \nu:T^*(\tilde X)\setminus \{0\}\rightarrow
\mathbb{C}$. The symplectic reductions of $\left (T^*(\tilde
X)\setminus \{0\}\right )\times \left (T^*(\tilde
  X)\setminus \{0\}\right )$
and $\left (T^*(X)\setminus \{0\}\right )\times \left
(T^*(X)\setminus
  \{0\}\right )$, given by
\begin{eqnarray*}
\left (\tilde \nu \times \tilde \nu \right )^{-1}(1,1)/
\mathbb{C}^*\times \mathbb{C}^*\, \supseteq \, (\nu \times \nu
)^{-1}(1,1)/S^1\times S^1\end{eqnarray*} are symplectically
equivalent to $T^*(\tilde M\times
  \tilde M)\setminus \{0\}$ and
$T^*(M\times M)\setminus \{0\},$ respectively. The $S^1$- and
$\mathbb{C}^*$-invariant functions $\zeta$ and $\tilde \zeta$
descend to functions $\zeta _r$ and $\tilde \zeta _r$ on
$T^*(M\times M)\setminus \{0\}$ and $T^*(\tilde M\times \tilde
M)\setminus \{0\}$, and $\tilde \zeta _r$ is the almost analytic
extension of $\zeta _r$. The reduced Lagrangian submanifold
$\Gamma _r\subseteq T^*(\tilde M\times \tilde M)\setminus \{0\}$
is defined as
$$\Gamma _r=:\left (\Gamma \cap (\tilde \nu \times \tilde
\nu )^{-1}(1,1)\right )/\mathbb{C}^*\times \mathbb{C}^*.$$ Let
$\Upsilon _r\subseteq T^*(\tilde M)\setminus \{0\}$ be the
reduction of $\Upsilon$, that is, the almost analytic submanifold
defined by $\tilde \zeta _r=0$. Given Lemma \ref{lem:unique},
$\Gamma _r$ may be charachterized geometrically as the unique
almost analytic Lagrangian submanifold of $T^*(\tilde M\times
\tilde M)\setminus \{0\}$ contained in $\Upsilon _r\times
\overline \Upsilon _r$ and with the given real locus, which is the
reduction of the real locus of $\Gamma$.

The pull-back of $\phi$ under $\chi \times \chi$ is the function
$$\phi '=i\left (1-\lambda \overline \mu \gamma (p,q) \right )$$
on $X\times X$, where $\gamma (p,q)=c(p,q)/\sqrt{\tilde \beta
(p)}\sqrt{\tilde \beta (q)}$. Hence we may view i), ii) and iii)
as conditions on $\gamma$.

Over $\mathrm{diag}(W)\subseteq M\times M\subseteq \tilde M\times
\tilde M$ the projection $\Gamma _r\rightarrow \tilde M\times
\tilde M$ is a local diffeomorphism. Therefore we can find,
locally along $\mathrm{diag}(W)$, a function $\tilde \gamma
(\tilde p,\tilde q)$ such that
$$\Gamma _r=\left \{\left (\tilde p, d_{\tilde p}\log (\tilde b),\tilde q,
-d_{\tilde q}\log (\tilde b)\right ) \, : \, \tilde p, \tilde q
\in \tilde M\right \}.$$ The real locus $(\Gamma _r)_\mathbb{R}$
of $\Gamma _r$ is the reduction under $S^1\times S^1$ of the real
locus of $\Gamma$; given Lemma \ref{lem:composition} and
(\ref{eqn:circ}), in local coordinates we have
\begin{eqnarray}\label{eqn:gammar}
(\Gamma _r)_{\mathbb{R}}= \left \{ \left ((p,\zeta ), h\cdot
(p,-\zeta)\right ): p\in W,\,\zeta \in T^*_p(M), \, \exists \,
\theta \mbox{ such that }\right . \nonumber \\
\left .\alpha _{(z,e^{i\theta})}=d\theta + \zeta , h\in H_{\Phi
(p)} \right \}.\end{eqnarray}

\begin{lem} \label{lem:hhorizontal}
If $x\in X$, $p=\pi (x)\in W$ and $h_s\in H_{\Phi (p)}$ is a
smooth path, the path $\varsigma  (s)=h_s\cdot x$ in $X$ is
horizontal.
\end{lem}

\noindent {\it Proof.} It suffices to show that if $p=\pi (x)\in
W$ and $\xi \in \frak{h}_f$, where $f=\Phi (p)$, then the induced
vector field $\xi _X$ on $X$ is horizontal at $x$, that is,
$\alpha _x \left (\xi _X(x)\right )=0$. In view of
(\ref{eqn:hor}), we are then reduced to proving that if $\xi \in
\frak{h}_f$ then $\phi _\xi (p)=0$. We have
$$\phi _\xi (p)=\left <\Phi (p),\xi \right >=<f,\xi>=0,$$
in view of (\ref{eqn:hf}).

\begin{lem} $\tilde \gamma$ is constant on $W\times _{p_W}W$.
Therefore, after adding a suitable additive constant to it, we may
assume that $\tilde \gamma =1$ on $W\times _{p_W}W$,  i.e. that
$c(p,q)=\sqrt{\beta (p)}\sqrt{\beta (q)}$ if $p\in W$ and
$q=h\cdot p$ with $h\in H_{\Phi (p)}$.
\label{lem:constant}\end{lem}

\noindent {\it Proof.} Let $(p_1,h_1\cdot p_1)\in W\times _{p_W}W$
be in the neighbourhood of $(p_0,p_0)$. Thus, $p_1\in W$ and
$h_1\in H_{\Phi (p_1)}$. We can join $(p_0,p_0)$ and
$(p_1,h_1\cdot p_1)$ within $W\times _{p_W}W$ going first from
$(p_0,p_0)$ to $(p_1,p_1)$ through a smooth path $\upsilon (t)=
(p_t,p_t)$ ($t\in [0,1]$) with $p_t\in W$ for every $t\in [0,1]$,
and then from $(p_1,p_1)$ to $(p_1,h_1\cdot p_1)$ through a path
$\nu(t)=(p_1,h_t\cdot p_1)$ ($t\in [0,1]$) with $h_0=e$ and
$h_t\in H_{\Phi (p_1)}$ for every $t\in [0,1]$. It is thus
sufficient to show that $\tilde b$ is constant on every path of
the above form, and for this in turn it suffices to check that
every element of $ (\Gamma _r)_{\mathbb{R}}$ vanishes on the
tangent vectors to these paths.

An element of $ (\Gamma _r)_{\mathbb{R}}$ lying over $(p_t,p_t)$
has the form $\eta =(p_t,\zeta, p_t,-\zeta)$ for a suitable
cotangent vector $\zeta$ at $p_t$, while $\left
.\frac{d\upsilon}{dt}\right |_t=(\dot p_t,\dot p_t)$. Thus, $\left
<\eta,(\dot p_t,\dot p_t)\right >=\zeta (\dot p_t)-\zeta (\dot
p_t)=0$.

An element of $ (\Gamma _r)_{\mathbb{R}}$ lying over
$(p_1,h_s\cdot p_1)$ has the form $$\eta =\left (p_1, \zeta ,
h_s\cdot p_1,-\zeta _s\right )$$ where $\zeta _s=\left
(d_{p_1}h_s^{-1}\right )^t(\zeta)$. Here by definition $\zeta \in
T^*_{p_1}M$ satisfies $\alpha _{(p_1,e^{i\theta})}= \zeta
+d\theta$, and $d_{p_1}h_s^{-1}$ denotes the differential of
$h_s^{-1}\in G$, viewed as a diffeomorphism of $M$,  at $p_1$.

On the other hand, $\left .\frac{d\nu}{dt}\right |_{t=s}= (0,\dot
\kappa (s))$, where $\kappa (s)=h_{s}\cdot p_1 $. Thus, we need to
check that $\zeta _s\left (\dot \kappa (s)\right )=0$ for every
$s$.

Let us consider, in the local trivialization induced by $\tilde
e_L^*$, the path in $X$ given for some fixed $\theta \in
\mathbb{R}$ by $$\tilde k(s)=h_s\cdot (p_1,e^{i\theta})=(h_s\cdot
p_1,e^{i\theta})$$ (the second equality reflects the construction
of $\tilde e_L$).

By Lemma \ref{lem:hhorizontal}, $\tilde k$ is horizontal, that is,
$ \alpha _{\tilde \kappa (s)} \left (\dot {\tilde \kappa}
(s)\right )=0$ for every $s$. Since, again invoking the definition
of $(\Gamma _r)_{\mathbb{R}}$ and the construction of $\tilde
e_L$, $\alpha _{\tilde \kappa (s)} =\zeta _s+d\theta$, we conclude
that $\zeta _s\left (\dot \kappa (s)\right )=0$.

\bigskip

Let now $\Gamma '\subseteq T^*\left (\tilde X\times \tilde X\right
)\setminus \{0\}$ be the almost analytic Lagrangian submanifold
parametrized by $\phi '$. By constructon, $\Gamma '$ satisfies
conditions (ii) and (iii) above. To verify condition (i), suppose
that $(x_0, y_0)\in Z(\phi ')$ and thus $\big (x_0, d_x\phi
'(x_0,y_0),y_0,-d_{y}\phi '(x_0,y_0)\big ) \in \Gamma
'_{\mathbb{R}}$ (the real locus of $\Gamma '$). Then (in view of
\cite{sz}, (44)) $d\log (\tilde \gamma )(x_0,y_0)$ must be real.
Thus, if $x_0,y_0\in X$ correspond to $\big ((p_0,\lambda
_0),(q_0,\mu _0)\big )$ under $\chi \times \chi$, we must have
$(p_0,q_0)\in W\times _{p_W}W$. Hence $\tilde \gamma (p_0,q_0)=1$,
and the condition that $\phi '(p_0,q_0)=1$ implies $\lambda _0=\mu
_0$. Thus, by our choice of local coordinates, if $h\in H_{\Phi
(p_0)}$ is such that $q_0=h\cdot p_0$, we have $y_0=h\cdot x_0$.
Hence, $\Gamma '$ satisfies (i) and the proof of Proposition
\ref{prop:phasegamma} is complete.

\begin{prop}\label{prop:positivetype}
The regular phase function of Proposition \ref{prop:phasegamma} is
of positive type, that is, the imaginary part of $\psi$ satisfies
\begin{eqnarray}
{\cal I}(\psi)\ge 0\label{eqn:positivetype}\end{eqnarray} on an
open neighbourhood of $(p_0,p_0)\in W\times _{p_W}W$.
\end{prop}

\noindent {\it Proof.} Let us go back to our original preferred
holomorphic coordinates in the neighbourhood of $p_0$. Recall the
geometric charachterization of $\Gamma$ described in Lemma
\ref{lem:unique}. The condition that $\Gamma \subseteq \Upsilon
\times \overline \Upsilon$ (Lemma \ref{lem:unique}, iii)) implies

\begin{claim} To infinite order along $W\times _{p_W}W$, $a$ is holomorphic in
$p$ and anti-holomorphic in
$q$.\label{claim:holantihol}\end{claim}

On the other hand, the condition that $\Gamma \supset \left
(D\circ
  \Sigma '\right )'$
be the real locus of $\psi$ (Lemma \ref{lem:unique}, ii)) forces
(in the given trivialization)

$$1=\lambda \overline \mu a(p,q)\,\,\Longleftrightarrow \,\,
\left ((p,\lambda),(q,\mu )\right ) \in D\circ \Sigma '.$$ In
particular, applying this to $\left ((p,\lambda),(p,\lambda
)\right )\in {\rm diag}(X)$ with $p\in W$, we obtain
$$a(p,p)\, =\, \beta (p)\,\,\,\,\,\,\,\,(p\in W).$$
In other words, $a(\cdot,\cdot)$ restricted to ${\rm diag}(W)$ is
the dual hermitian metric $\beta$ on $L^*$, in the holomorphic
coframe $e_L^*$.

In the local holomorphic coordinates $(z_i,w_j)$, with
$w_j=a_j+ib_j$, $W$ is defined by the equations $b_j=0$. Thus,
working in local coordinates near $p\in W$ and given
(\ref{eqn:expansion}), $a\left ((z,a+i0),(z,a+i0)\right
)=1+|z|^2+|a|^2+\cdots$. In view of Claim \ref{claim:holantihol}
we conclude

\begin{cor}
To infinite order along $W$,
$$a\left ((z,w),(z',w')\right )=1+z\cdot \overline {z'}+w\cdot
\overline {w'}+\cdots.$$ \label{cor:expansion}\end{cor}

We can now apply {\it verbatim} the argument following the proof
of Lemma 2.2 in \cite{sz} (we only need to replace $z$ with the
pair $(z,w)$) to conclude that $\psi$ is a positive phase
function. More precisely, let us fix $x_0\in \pi ^{-1}(p_0)$. We
may assume that $\tilde e_L^*(p_0)=x_0$, so that $x_0$ has local
coordinates $(z_0,w_o,\theta _0)=(0,0,0)$, where $\lambda
=e^{i\theta}$). A nearby $x\in X$ has local coordinates
$(z,w,\theta)$. Let us compute ${\cal I}(\psi (x_0,x))$ as the
real part of $-i\psi (x_0,x)$. We obtain
\begin{eqnarray*}
-i\, \psi (x_0,x) \, =\, 1-\frac{a(0,0,z,w)}{\sqrt {\beta
(z,w)}}e^{-i\theta}
=\left (1-e^{-i\theta}\right )\\
+e^{-i\theta}\left [\frac 12 \|(z,w)\| ^2+ O\left (
\|(z,w)\|^3\right )\right ],\end{eqnarray*} which is positive if
$\theta$ and $(z,w)$ are small.

\bigskip

Given this, we may conclude:

\begin{prop}
$P_{(\omega)}$ may be microlocally represented near the diagonal
by an oscillatory integral of the kind:
\begin{eqnarray*}
\int _0^\infty \,e^{it\psi _\omega (x,y)}\, s_\omega
(x,y,t)dt,\end{eqnarray*} where $s_\omega \in S^n(X\times X\times
\mathbb{R}^+)$ has an asymptotic development
$$s_\omega (x,y,t)\sim \sum _{j=0}^\infty t^{n-j}\,s_{\omega}^{(j)}(x,y).$$
By ellipticity $s_{\omega}^{(0)}(x,y)\neq 0$ if $(x,y)\in \Gamma
_\mathbb{R}$.
\end{prop}

With this microlocal description, to complete the proof of Theorem
\ref{thm:2} we need now only follow the argument in the proof of
Thoerem 1 on pp 327-9 of \cite{z}.

Corollaries \ref{cor:1} and \ref{cor:2} are now straightforward,
since the complexification $\tilde G$ of $G$ acts on each
irreducible piece $H^0(M,L^{\otimes k})_\omega$.

\bigskip

Let us now remove the simplifying assumption that $\omega$ lie in
an elementary fundamental wedge for the Hamiltonian action of $G$
on $T^*X$. Recall that $\Psi :T^*X\rightarrow \frak{g}^*$ is the
moment map.

\begin{lem} Let $R\subseteq T^*X\setminus \{0\}$ be a $G$-invariant
open conic set satisfying the following conditions:

\noindent i): the restriction $\left .\Psi \right |_R:R\rightarrow
\frak{g}^*$ is trasversal to $C(\mathcal{O})$;

\noindent ii): $G_\omega$ acts freely on $\Psi ^{-1}(\omega)\cap
R$.

Then $Z_R=\Psi ^{-1}\big (C(\mathcal{O})\big )\cap R$ is a
fibrating coisotropic submanifold of $R$, and the leaf of the null
fibration through $y\in Z_R$ is the orbit of $H_f$, $f=\Psi (y)$.
Let $p_{Z_R}:Z_R\rightarrow Z_R^\sharp$ be the null fibration.
Then the fibre product $Z_R\times _{p_{Z_R}}Z_R\subseteq
T^*X\times T^*X$ is a Lagrangian relation. Furthermore, the
orthogonal projector $Q_{(\omega)}:L^2(X)\rightarrow
H(X)_{(\omega)}$ is microlocally equivalent on $R\times R\subseteq
T^*X\times T^*X$ to a Fourier integral operator associated to
$Z_R\times _{p_{Z_R}}Z_R$. \label{lem:micro}
\end{lem}

In other words, there is a Fourier integral operator $\tilde
Q_{(\omega)}$ on $X$ associated to the Lagrangian relation
$Z_R\times _{p_{Z_R}}Z_R$ such that $$\mathrm{WF}\left
(Q_{(\omega)}- \tilde Q_{(\omega)}\right )\cap (R\times
R)=\emptyset.$$

\noindent \textit{Proof.} This follows straightforwardly from the
theory of \cite{gs-hq}, but we shall sketch a proof for the
reader's convenience.

The first part of the Lemma is an immediate consequence of the
arguments in section 2 of \textit{loc. cit.}.

On the other hand,
\begin{equation}\label{eqn:momlag}
\Gamma =\left \{(g,\gamma,\zeta,g\zeta)\,:\,g\in G,\,\zeta\in
T^*X\setminus \{0\},\, \gamma =\Psi (\zeta)\right \}\end{equation}
is a canonical relation in $T^*G\times T^*X\times (T^*X)^-$ (cfr
(4.2) of \textit{loc. cit.}), called the \textit{moment
Lagrangian} for the action of $G$ on $T^*X$. In our case, this is
simply the conormal bundle to the graph of the action of $G$ on
$X$. Here we implicitly identify $T^*G$ with $G\times \frak{g}^*$
by means of right translations, and denote by $(T^*X)^-$ the
symplectic manifold obtained by $T^*X$ by changing sign to the
standard symplectic structure. Restricting to $R\times R$, we
obtain a homogeneous canonical relation
\begin{equation}\label{eqn:momlagR}
\Gamma _R=\left \{(g,\gamma,\zeta,g\zeta)\,:\,g\in G,\,\zeta\in
R,\, \gamma =\Psi (\zeta)\right \}\subseteq T^*G\times R\times
R^-.
\end{equation}

Let furthermore $$\Lambda _{C(\mathcal{O})}=\{(g,f):f\in
C(\mathcal{O}), g\in H_f\}\subseteq G\times \frak{g}^*\cong T^*G$$
be the charachter Lagrangian associated to $C(\mathcal{O})$, a
homogeneous Lagrangian submanifold of $T^*G$. Under the hypothesis
of the Lemma, $\Gamma _R$ and $\Lambda _{C(\mathcal{O})}$
intersect cleanly, and $Z_R\times _{p_{Z_R}}Z_R=\Gamma _R^t\circ
\Lambda _{C(\mathcal{O})}$.

Now recall that $$Q_{(\omega)}=\int _G\chi _{(\omega)}
(g^{-1})\varrho (g)\,dg,$$ where by $\varrho$ we denote the
unitary action of $G$ on $L^2(X)$, and $$\chi _{(\omega)} =\sum
_{k=1}^\infty \dim (V_{k\omega}) \chi _{k\omega}.$$ Here $\chi
_{k\omega}$ is the charachter function associated to the
irreducible representation with highest weight $k\omega$. By
Theorem 6.3 of \cite{gs-hq}, $\chi _{(\omega)}$ is a Lagrangian
distribution on $G$, associated to the charachter Lagrangian
$\Lambda _{C(\mathcal{O})}$; on the other hand, the Schwartz
kernel $\tilde \varrho \in \mathcal{D}' \left (G\times X\times
X\right )$ of the action $\varrho$ is a Lagrangian distribution on
$T^*G\times T^*X \times (T^*X)^-$, associated to the moment
Lagrangian. By the above, the Schwartz kernel $\tilde
Q_{(\omega)}\in \mathcal{D}'(X\times X)$ of $Q_{(\omega)}$ is
given by:
$$\tilde Q_{(\omega)}(x,y)=\int _G\chi _{(\omega)}
(g^{-1})\tilde \varrho (g,x,y)\,dg.$$ Let $P_1+P_2\sim
\mathrm{id}$ be a pseudodifferential partition of unity on $X$
such that $P_1\sim \mathrm{id}$ on $R$. We may write $\tilde
\varrho =\tilde \varrho _1+\tilde \varrho _2$, where $\tilde
\varrho _j=P_j \circ \tilde \varrho $, and accordingly $\tilde
Q_{(\omega)}\sim \tilde Q_{(\omega),1}+\tilde Q_{(\omega),2}$.
Then $\tilde Q_{(\omega),2}$ is smoothing on $R\times R$. Since on
the other hand $\tilde \varrho _1$ is a Lagrangian distribution
associated to $\Gamma _R$, by the hypothesis we may interpret
$\tilde Q_{(\omega),1}$ as a pull-back (under the diagonal map
$G\times X\times X\rightarrow G\times G\times X\times X$),
followed by a push-forward. The Lemma then follows from the usual
functorial properties of distributions.

\bigskip

Let us now return to the proof of the Theorem. Given the relation
between $\Psi$ and $\Phi$ on the cone $Y\subseteq T^*(X)$, if
$\Phi$ is transversal to $C(\mathcal{O})$ then at any rate so is
$\left .\Psi \right |_Y$, whence $\Psi$ itself in a conic open
neighbourhood $R$ of $Y$ in $T^*X\setminus \{0\}$. Given the
$G$-equivariance of $\Psi$, we may assume that $R$ is
$G$-invariant. By assumption, $G_\omega$ acts freely on $\Phi
^{-1}\left (\mathbb{R}_+\omega\right )$, whence on $\Psi
^{-1}\left (\mathbb{R}_+\omega\right )\cap Y$. Therefore, perhaps
after replacing $R$ by a smaller $G$-invariant conic open
neighbourhood of $Y$ in $T^*X\setminus \{0\}$, we may assume that
$G_f$ acts freely on $\Psi ^{-1}\left (\mathbb{R}_+\omega\right
)\cap R$. Arguing as above with a pseudodifferential partition of
unity, we may decompose $\tilde Q_{(\omega)}$ (or $Q'_{(\omega)}$)
as a sum $\tilde Q_{(\omega)}=\tilde Q_{(\omega),1}+\tilde
Q_{(\omega),2}$, where $\tilde Q_{(\omega),2}$ is smoothing on
$R\times R$. Therefore, so is the composition $Q_{(\omega),2}\circ
\Pi$, since $\mathrm{WF}(\Pi)\subseteq Y\times Y\subseteq R\times
R$. Every choice in the above argument may be assumed
$S^1$-invariant. We are thus reducing to applying the previous
arguments to the composition $Q_{(\omega),1}\circ \Pi$.

\section{Proof of Proposition \ref{prop:rapiddecay}.}

\noindent The proof is a slight modification of some of the
arguments in the proofs of Theorem 1 and of Theorem \ref{thm:3}.
Fix $p\in M$ with $\Phi (p)\neq 0$, and for $\epsilon >0$ let
$$U_{p,\epsilon}=:\{g\in G: {\rm dist}_M(gp,p)<\epsilon \}.$$
Then $U_{p,\epsilon }$ is an open neighbourhood of the stabilizer
subgroup $$G_p=\{g\in G:gp=p\},$$ and $U_{p,\epsilon }\rightarrow
G_p$ as $\epsilon \rightarrow 0$. We have $G_p=\{e\}$ if $p\not
\in R$, since in our situation $G$ acts freely on $M\setminus
R=\tilde G\cdot \Phi ^{-1}(0)$.

Choose $x\in X$ with $p=\pi (x)$. Then $\Pi _{k,\omega}(p,p)=\Pi
_{k,\omega}(x,x)$ may be decomposed as in (\ref{eqn:komega1}) with
$U=U_{p,\epsilon }$; thus the first term is the integral of an
averaged Szeg\"{o} kernel over $U_{p,\epsilon }$ and the second is
the integral over $G\setminus U_{p,\epsilon }$. For any fixed
$\epsilon $, the latter term is therefore rapidly decreasing as
$k\rightarrow +\infty$, since $\Pi _k(x',x^{\prime\prime})$ is
rapidly decreasing in $k$ if $\pi (x')\neq \pi (x^{\prime
\prime})$. We are then led to estimate
\begin{eqnarray}\Pi _{k,\omega}(x,x)'=:
\dim (V_\omega )\cdot \int _{U_{p,\epsilon }} \, \chi _\sigma
(g^{-1})\,\Pi _k(\tilde \mu _{g^{-1}}(x),x)\, dg,
\end{eqnarray}
and we shall now argue that if $\epsilon$ has been chosen suitably
small (in particular, very small with respect to $\varphi (p)
=\left \|\Phi (p)\right \|>0$) then $\Pi _{k,\omega}(x,x)'$ is
also rapidly decreasing as $k\rightarrow +\infty$. To this end,
let us first remark that, as in (\ref{eqn:komega2}), $\Pi
_{k,\omega}(x,x)'$ may be written as a complex oscillatory
integral in $dg,\,d\theta,\, du$, with phase $\Psi _x(g,\theta,u)
\, = \, u \psi (e^{i\theta}\mu _{g^{-1}}(x),x) - \theta$ (here and
below we think of $x$ as a parameter, and denote by $d\Psi _x$ the
differential of $\Psi _x$ with respect to $(g,\theta,u)$).

Arguing as in the first part of the proof of Lemma \ref{lem:psix},
we can see that given $0<a<1<b$ there exists $d>0$ and $\epsilon
'>0$ such that if $0<\epsilon <\epsilon '$ then $\|d\Psi _x
(g,\theta,u)\|>d$ if $g\in U_{p,\epsilon }$ and $u\not \in [a,b]$.

Now if $g\in U_{p,\epsilon }$ we have an estimate similar to
(\ref{eqn:estimateondpsi}), with $\epsilon _0^2$ replaced by
$\epsilon $, and to (\ref{eqn:estimateonalfa}), with $C_2\epsilon
_0$ replaced by $\varphi (p)$. Therefore, if $\epsilon $ has been
chosen sufficiently small, the phase $\Psi _x$ has no stationary
points on $U_{x,\epsilon }\times S^1\times (0,+\infty)$, and
actually its differential is bounded from below in norm there.
This proves the first statement of the proposition.

Now let $K\subset M\setminus \left (\Phi ^{-1}(0)\cup R\right )$
be a compact subset, and let us show that the rapid decay of $\Pi
_{k,\omega}(x,x)'$ is uniform over $x\in \pi ^{-1}(K)$. For
$\delta
>0$ let $U_\delta \subseteq G$ be the $\delta$-neighbourhood of
the identity in $G$, with respect to some fixed invariant metric.
As $K$ is compact and $G$ acts freely on the semistable locus
$\tilde G\cdot \Phi ^{-1}(0)\supset K$, there exist
$0<r_{K,\delta}<s_{K,\delta}$ such that
\begin{equation}\mathrm{dist}_X(gx,x)<s_{K,\delta}\text{ for any
}x\in X \text{ and } g\in
U_\delta,\label{eqn:skdelta}\end{equation}
\begin{equation}\mathrm{dist}_X(gx,x)>r_{K,\delta}\text{ if }\pi (x)\in K
\text{ and } g\not\in U_\delta.\label{eqn:rkdelta}\end{equation}
Furthermore, we may assume that
$r_{K,\delta}<s_{K,\delta}\rightarrow 0\text{ as }\delta
\rightarrow 0$.

Let us then choose $\delta >0$ so that $s_{K,\delta}$ is very
small with respect to $$\varphi (K)=:\min \{\varphi (p):p\in
K\}>0.$$ For every $x\in \pi ^{-1}(K)$ we may decompose $\Pi
_{k,\omega}(x,x)$ in (\ref{eqn:komega1}), where the first term is
now the integral over $U_\delta$ and the second is the integral
over $G\setminus U_\delta$. In view of (\ref{eqn:rkdelta}), the
latter term is rapidly decreasing in $k$ as $k\rightarrow
+\infty$.

We are thus led to look for a uniform estimate, over $x\in \pi
^{-1}(K)$, of
\begin{eqnarray}\Pi _{k,\omega}(x,x)'=:\dim (V_\omega )\cdot \int _{U_\delta} \,
\chi _\sigma (g^{-1})\,\Pi _k(\tilde \mu _{g^{-1}}(x),x)\, dg.
\end{eqnarray}
Viewing this as a complex oscillatory integral as in
(\ref{eqn:komega2}), we want to bound from below the differential
of the phase $\Psi _x(g,\theta,u) \, = \, u \psi (e^{i\theta}\mu
_{g^{-1}}(x),x) - \theta$, uniformly over $x\in \pi ^{-1}(K)$.

Given that $\pi ^{-1}(K)$ is compact, a slight modification of the
argument used in the proof of Lemma \ref{lem:psix} now proves the
following: If $0<a<1<b$ and $0<c<\pi$ there exists $\gamma
>0$ such that for all sufficiently small $\delta >0$ we have
$\|d\Psi _x\|>\gamma$ if $x\in \pi ^{-1}(K)$, $g\in U_\delta$ and
$u\not \in [a,b]$ \textit{or} $\theta \not \in [-c,c]$.

We may now argue as in the first part of the proof of this
Proposition, however replacing now $\epsilon _0^2$ with
$s_{K,\delta}$ in the analogue of (\ref{eqn:estimateondpsi}),
holding if $(x,g)\in K\times U_\delta$, and $C_2\epsilon _0$ with
$\varphi (K)$ in the analogue of (\ref{eqn:estimateonalfa}).

\end{document}